\documentclass{gtpart}     
\agtart

\usepackage{amssymb, amsmath, amsfonts}
\usepackage{xypic}
\xyoption{all}

\title{Representations of Spaces}

\author{Wojciech Chach\'olski}
\givenname{Wojciech} \surname{Chach\'olski}
\address{K.T.H.\\ Matematik\\ S-10044 Stockholm,
Sweden} \email{wojtek@math.kth.se} \urladdr{}

\author{J\'er\^ome Scherer}
\givenname{J\'er\^ome} \surname{Scherer}
\address{Universitat Aut\`onoma de Barcelona \\ Departament de Matem\`atiques\\
E-08193 Bellaterra, Spain} \email{jscherer@mat.uab.es} \urladdr{}



\arxivreference{math.AT/0511577} \arxivpassword{}

\volumenumber{} \issuenumber{} \publicationyear{} \papernumber{}
\startpage{}
\endpage{}
\doi{} \MR{} \Zbl{} \received{} \revised{} \accepted{}
\published{} \publishedonline{} \proposed{} \seconded{}
\corresponding{} \editor{} \version{}

\newtheorem{theorem}{Theorem}[section]
\newtheorem{lemma}[theorem]{Lemma}
\newtheorem{proposition}[theorem]{Proposition}
\newtheorem{corollary}[theorem]{Corollary}
\newtheorem{point}[theorem]{\hspace{-1.7mm}}

\newtheorem{definition}[theorem]{Definition}
\newtheorem{example}[theorem]{Example}

\newtheorem{remark}[theorem]{Remark}

\def\mathbold#1{#1} 

\newcommand{\colim}[2]{\mbox{$\text{\rm colim}_{#1}#2$}}

\def\C{{\mathcal {C}}}
\def\D{{\mathcal {D}}}
\def\M{{\mathcal {M}}}

\def\rmono{\rto|<\hole|<<\ahook}
\def\rrmono{\rrto|<\hole|<<\ahook}
\def\umono{\ar@{_{(}->}[u]}
\def\uumono{\ar@{_{(}->}[uu]}
\def\dmono{\dto|<\hole|<<\ahook}
\def\drmono{\drto|<\hole|<<\ahook}

\def\lmono{\ar@{_{(}->}[l]}
\def\llmono{\ar@{_{(}->}[ll]}

\def\lra{\longrightarrow}

\def\ra{\rightarrow}
\def\la{\leftarrow}
\def\xra{\xrightarrow}
\def\xla{\xleftarrow}

\def\mono{\hookrightarrow}

\begin{document}

\begin{abstract}
We explain how the notion of homotopy colimits gives rise to that
of mapping spaces, even in categories which are not simplicial. We
apply the technique of model approximations and use elementary
properties of the category of spaces to be able to construct
resolutions. We prove that the homotopy category of any monoidal
model category is always a central algebra over the homotopy
category of $Spaces$.
\end{abstract}

\maketitle


\section{Introduction}\label{sec intro}

Let $\C$ be a category with a chosen class of weak equivalences
$W$ ($W$ is a collection of morphisms in $\C$ that satisfies the
``2 out of 3'' property and contains all isomorphisms). What does
it mean that one can do homotopy theory on $(\C,W)$? For us it
means three things. First, we should be able to perform basic
operations in $\C$ such as push-outs/pull-backs and more generally
colimits/limits. Second, the derived versions of these functors
(hocolims and holims) should exist. Third, these derived functors
should yield a certain action of the category of simplicial sets
$Spaces$ on $(\C,W)$ and lead to mapping spaces. Our attitude  is
that colimits/limits and their derived functors are the
fundamental ingredients of homotopy theory and the rest should
follow from their properties. This is an alternative  to the
approach of Hovey presented in~\cite{99h:55031} and more in the
spirit of what we believe we have learned from Bousfield, Dwyer,
Kan and others, \cite{MR51:1825} and \cite{wgdkanhirsmi}.

In this paper we study the following functor:
\[
\otimes_{l}:Spaces\times \C\ra \text{Ho}(\C)\ \ \ \ {\mathbf
K}\otimes_{l}X:= \text{hocolim}_{\mathbf K}cX
\]
where $\text{Ho}(\C)=\C[W^{-1}]$ and $cX:\mathbf K\ra \C$ is the
constant functor indexed by the simplex category of $K$
(see~\ref{pt simplexcat}) whose value is $X$. We would like
$\otimes_{l}$ to be the composition of the localization functor
$Spaces\times \C\ra\text{Ho}(Spaces)\times \text{Ho}(\C)$ and a
functor denoted by the same symbol
$\otimes_{l}:\text{Ho}(Spaces)\times \text{Ho}(\C)\ra
\text{Ho}(\C)$ (homotopy invariance). The latter should be a
homotopy left action (see~\ref{pt leftsimpact} and~\ref{pt
hosimpstruct}), which implies in particular that, for any $X\in
\text{Ho}(\C)$, the functor  $-\otimes_{l} X:\text{Ho}(Spaces)\ra
\text{Ho}(\C)$ has a right adjoint $\text{map}(X,-)$. This is in
our view the universal characterization of what mapping spaces
should be.

How can this be realized? In order to insure the existence of the
localized category $\C[W^{-1}]$ and the derived functors
$\text{hocolim}$s and $\text{holim}$s satisfying the above
conditions, we need to be able to do some homotopical algebra on
$(\C,W)$. One way of formalizing this is to assume that $(\C,W)$
can be given a Quillen (simplicial) model
structure~\cite{MR36:6480}. Dwyer and Kan explained
in~\cite{MR81m:55018} how to enrich such a  category with an
action of $Spaces$ using the so called ``hammocks''. They also
introduce the technique of taking (co)simplicial resolution of the
source (target), further developed by Hovey in \cite{99h:55031},
as well as Hirschhorn in \cite{MR1944041} to obtain mapping
spaces. Alternatively one could try to find conditions under which
a model category is Quillen equivalent to a simplicial model
category, as defined by Quillen in \cite[Section~II.1]{MR36:6480}.
This has been done by Dugger in \cite{MR2002f:55043} and Rezk,
Shipley, and Schwede in \cite{MR2002d:55025}. These approaches
however do not seem to emphasize the fundamental role of homotopy
colimits and limits. In this context we should mention also the
work of Heller~\cite{MR920963}.

Since in our view hocolims and holims are so important, we need a
suitable set up in which such constructions can be effectively
studied. It turns out that for that purpose, instead of putting
additional assumptions on a model category
(cellularity~\cite{MR1944041}, realization
axiom~\cite{MR2002d:55025}, etc...), it is more advantageous to
relax some of the requirements. The aim of this paper is to
explain the less restrictive approach of model approximations as
introduced in~\cite{MR2002k:55026}. Recall:

\medskip

\begin{definition}\label{def modelapprox}
{\em A {\it left model approximation} of $(\C,W)$ is a pair of
adjoint functors  $l:\M \rightleftarrows \C:r$ satisfying the
following conditions:
 \begin{enumerate}
\item $\M$ is a model category;
\item $l$ is left adjoint to $r$;
\item if $f$ is a weak equivalence in $\C$, then  $rf$ is a weak equivalence in $\M$;
\item if $f$ is a weak equivalence between cofibrant objects in $\M$, then
$lf$ is a weak equivalence in $\C$;
\item for any $A\in \C$ and cofibrant $X\in \M$, if $X\ra rA$ is a weak equivalence in $\M$,
then its left adjoint $lX\ra A$ is a weak equivalence in $\C$.
\end{enumerate}

A pair of adjoint functors $l:\C \rightleftarrows \M:r$ is a right
model approximation of $(\C,W)$ if the duals
$r^{\vee}:\M^{op}\rightleftarrows \C^{op} :l^{\vee}$ form a left
model approximation of $(\C^{op},W)$.}
\end{definition}

Model approximations have several crucial properties. A category
admitting a model approximation (left or right) can be localized
with respect to weak equivalences to form the homotopy category
$\text{Ho}(\C)$~\cite[Proposition 5.5]{MR2002k:55026}. Weak
equivalences in such categories are saturated, i.e., a morphism in
$\C$ belongs to $W$ if and only if the induced morphism in
$\text{Ho}(\C)$ is an isomorphism. The category $\text{Ho}(\C)$
can be identified with a full subcategory of $\text{Ho}(\M)$ via
the functor $r$ or $l$, depending on whether the approximation is
left or right (\cite[Proposition~5.5]{MR2002k:55026}). A left
model approximation $l:\M \rightleftarrows \C:r$ leads to a
``cofibrant replacement'' in $\C$ given by the adjoint $lQrA\ra A$
of a cofibrant replacement $QrA\ra rA$ in
$\M$~\cite[Remark~5.10]{MR2002k:55026}. In a dual way a right
model approximation leads to a fibrant replacement. These
replacements can be used to construct derived functors. If
$F:\C\ra \D$ is a functor for which the composition $Fl$ sends
weak equivalences between cofibrant objects in $\M$ to
isomorphisms in $\D$, then the left derived functor of $F$ is
given by $A\mapsto F(lQrA)$~\cite[Proposition 5.9, Remark
5.10]{MR2002k:55026}. In the same way right model approximations
can be used to construct right derived functors. A left model
approximation $l:\M \rightleftarrows \C:r$ induces a natural left
model approximation on the category of functors
$\text{Fun}(I,\C)$, for any small category $I$ (see
Theorem~\ref{thm hocolimconst}). A cofibrant replacement of this
model approximation is thus good for constructing the total left
derived functors of the colimit functor
$\text{colim}_{I}:\text{Fun}(I,\C)\ra \C$. Likewise one constructs
the total left derived functor of left Kan extensions,
\cite[Theorem~11.3(3)]{MR2002k:55026}. These derived functors (the
homotopy colimits and the homotopy left Kan extensions) are
related to each other by various properties. At a $2$-categorical
level this has been formalized by Cisinski in~\cite{Cisinski1} and
amounts to saying that the prederivator associated to any model
category is a Grothendieck derivator. The same applies to right
model approximations.

To make the approach of model approximations viable we need one
more thing. We need to show that a model approximation leads in a
natural way to an action of $Spaces$. The purpose
of~\cite{MR2002k:55026} was to introduce simple geometric
techniques to study hocolims and holims. In this paper we are
going to illustrate how to use these techniques further, which is
what motivated us originally.  We did not want to just state and
prove some properties of mapping spaces, but more importantly we
meant to explain an approach and provide coherent tools to study
them. This will be illustrated by proving the following theorem
(see Proposition~\ref{prop mapasadjoint} and ~\ref{prop
homotopyadjoint}):

\medskip

{\bf Theorem.}
{\it Assume that $(\C,W)$ admits a left model approximation. The
functor $\otimes_{l}:Spaces\times \C\ra \text{\em Ho}(\C)$ is
homotopy invariant and, for any $X\in \C$,
$-{\otimes_{l}}X:\text{\em Ho}(Spaces)\ra \text{\em Ho}(\C)$ has a
right adjoint $\text{\em map}(X,-)$.

Assume that $(\C,W)$ admits both a left and a right model
approximations. Then, for any space $K$, the functor $\mathbf
K\otimes_{l}-:\text{\rm Ho}(\C)\ra \text{\rm Ho}(\C)$ is left
adjoint to $X\mapsto \text{\rm holim}_{\mathbf K}cX$.}

\medskip

As an application we prove a conjecture made by Hovey in
\cite[p.~119]{99h:55031} about the homotopy category of a monoidal
model category. This has been proved independently by Cisinski in
his work on Grothendieck derivators, \cite{Cisinski}.



\medskip

The plan of the article is as follows. Section~\ref{sec notation}
is a list of basic definitions and notation. In particular we
recall what an action of $Spaces$ on a model category is, and how
to construct derived functors. The reader might skip this section
and refer to it when needed. In Section~\ref{sec simplexcat} we
explain that, unlike arbitrary categories, simplex categories
behave well with respect to pull-backs. A brief reminder about
bounded functors appears then in Section~\ref{sec bounded}. We
introduce in Section~\ref{sec subdivision} a simplicial structure
on bounded diagrams indexed by a barycentric subdivision. In
Section~\ref{sec hobounded} we see how the geometry of the simplex
category is reflected in the homotopy theory of bounded diagrams.
The next three sections form the core of the paper.
Section~\ref{sec constatfun} is basically devoted to prove an
invariance property for homotopy colimits of constant diagrams. In
Section~\ref{sec almost} we show that very little is missing for
the simplicial structure to be homotopy invariant; for
homotopically constant diagrams, it is so. We are finally ready to
construct mapping spaces in Section~\ref{sec resolution}. The
theory is dualized in Section~\ref{sec right}. In the final
Section~\ref{sec monoidal} we illustrate our approach by giving an
answer to Hovey's question.

\noindent {\bf Acknowledgments.} The first author is supported in
part by NSF grant DMS-0296117, Vetenskapsr\aa det grant 2001-4296,
and G\"oran Gustafssons Stiftelse, the second author by FEDER-MEC
grant MTM2004-06686 and the program Ram\'on y Cajal, MEC, Spain.

\section{Notation and set up}
\label{sec notation}

\begin{point}\label{pt cat}{\rm
$Cat$ denotes the category of small categories and $Spaces$ the
category of simplicial sets. A simplicial set is also called a
space.}
\end{point}

\begin{point}\label{pt limcolim}{\rm
From Section 3 on, categories denoted by the symbols $\C$ and $\D$
are assumed to contain all colimits and limits, compare with
Quillen's axiom (MC0) \cite{MR36:6480}.}\end{point}

\begin{point}{\rm
A category with weak equivalences is a category $\C$ with a chosen
collection  of morphisms $W$, called weak equivalences, which
satisfies the ``2 out of 3'' property and contains all
isomorphisms. Elements in $W$ are denoted by the symbol
$\xra{\sim}$. The category $\C[W^{-1}]$, if it exists, is called
the homotopy category  and is denoted by $\text{\rm Ho}(\C)$. By
$\pi_{\C}:\C\ra \text{\rm Ho}(\C)$ we denote the localization
functor. }\end{point}

\begin{point}\label{pt functors}{\rm
The category of all $\C$-valued functors indexed by a small
category $I$, with natural transformations as morphisms, is
denoted by $\text{Fun}(I,\C)$. The symbol $cX:I\ra \C$ denotes the
constant functor indexed by $I$ with value $X\in \C$. If $(\C,W)$
is a category with weak equivalences, then $\text{Fun}(I,\C)$
stands for a category with weak equivalences given by the
object-wise weak equivalences. A functor $F:I\ra \C$ is called
homotopically constant if there is an object $X\in\C$ and a weak
equivalence $F\ra cX$ in $\text{Fun}(I,\C)$. }\end{point}

\begin{point}{\rm
Let $I$ and $J$ be small categories and $f:I\ra J$ be a functor.
The left and right adjoints to $f^{\ast}:\text{Fun}(J,\D)\ra
\text{Fun}(I,\D)$, $F\mapsto F \circ f$, if they exist, are
 called  Kan extensions. They are denoted
respectively by $f^{k}: \text{Fun}(I,\D)\ra \text{Fun}(J,\D)$ and
$f_{k}: \text{Fun}(I,\D)\ra \text{Fun}(J,\D)$. These functors can
be described explicitly in terms of the over-categories
$f\downarrow j$ and under-categories $j\downarrow f$
~\cite[II.6]{MR2001j:18001}:
$(f^{k}F)(j)=\text{colim}_{f\downarrow j}F$ and
$(f_{k}F)(j)=\text{lim}_{j\downarrow f}F$. }\end{point}

\begin{point}\label{pt simplexcat}{\rm
The simplex category of a simplicial set $K$, denoted by the same
letter in a bold-face font~$\mathbf K$, is a category whose
objects are simplices of $K$, or equivalently maps of the form
$\sigma:\Delta[n]\ra K$. The morphisms in $\mathbf K$ between
$\sigma:\Delta[n]\ra K$ and $\tau:\Delta[m]\ra K$ are those maps
$f:\Delta[n]\ra \Delta[m]$ for which $\sigma=\tau f$. Subject to
standard cosimplicial identities, the morphisms in $\mathbf K$ are
generated by the face and degeneracy morphisms:
\[\xymatrix{\Delta[n]\rrto^{d_{i}} \drto & & \Delta[n+1]\dlto\\
& K }\ \ \ \ \
\xymatrix{\Delta[n+1]\rrto^{s_{i}} \drto & & \Delta[n]\dlto\\
& K }
\]
A map of spaces $f:A\ra B$ induces a functor $\mathbf f:\mathbf A
\ra \mathbf B$. For any simplex $\sigma \in B$, the functor
${\mathbf f}\downarrow \sigma\ra \mathbf A$ is a functor between
simplex categories given by the map of spaces
$\text{lim}(\Delta[n]\stackrel{\sigma}{\ra}B\stackrel{f}{\la}A)\ra
A$. In particular, when $f$ is the identity, the functor ${\mathbf
B}\downarrow \sigma\ra \mathbf B$ is induced by  the map
$\sigma:\Delta[n]\ra B$ }\end{point}

\begin{point}\label{pt simplexcatcolim}{\rm
The functor $Spaces\ra Cat$,  $K\mapsto \mathbf K$,  has a right
adjoint and therefore converts colimits in $Spaces$ into colimits
in $Cat$, see~\cite{MR516607}. This functor also takes pull-backs
of spaces into pull-backs of categories. It does not however
preserve final objects and consequently neither products:  the
simplex category of a point is the category of finite ordered sets
$\Delta$. For that reason we use the symbol $\mathbf A\boxtimes
\mathbf B$ to denote the simplex category of $A\times B$. The
product of the simplex categories is denoted as usual by ${\mathbf
A}\times {\mathbf B}$. The natural inclusion $\mathbf A\boxtimes
\mathbf B\ra \mathbf A\times \mathbf B$ is never an equivalence of
categories, unless $A$ or $B$ is the empty space. }\end{point}

\begin{point}\label{pt nerve}{\rm
$N:Cat\ra Spaces$ denotes the nerve
construction~\cite[Definition~6.3]{MR2002k:55026}. It has a left
adjoint and thus it converts limits in $Cat$ into limits in
$Spaces$, see~\cite{MR516607}. }\end{point}

\begin{point}\label{pt subdiv}{\rm
The functor $Spaces\ra Spaces $, $K\mapsto N({\mathbf K})$, is
called the barycentric subdivision. It  commutes with colimits
since it is a left adjoint to the functor that associates to a
space $L$, the space whose $n$-simplices are $mor_{Spaces}
(N({\mathbf \Delta}[n]), L)$, \cite[6.11]{MR2002k:55026}, which is
very similar to the so called Kan's $Ex$ construction,
\cite{MR19:759e}. }\end{point}

\begin{point}\label{pt wecat}{\rm
A functor $f:I\ra J$ between small categories is called a weak
equivalence if $N(f)$ is a weak equivalence in $Spaces$. A map of
spaces $f:A\ra B$ is a weak equivalence if and only if the induced
functor of simplex categories is a weak equivalence. For any
spaces $A$ and $B$, the inclusion ${\mathbf A}\boxtimes {\mathbf
B}\ra {\mathbf A}\times {\mathbf B}$ is a weak equivalence of
categories. }\end{point}

\begin{point}\label{pt leftsimpact}{\rm
A left action of $Spaces$ on $\C$ is an
 enrichment of $\C$ over
$Spaces$ for which the space of morphisms has a left adjoint in
the first variable. A left action of $Spaces$ on $\C$ can be
described explicitly as a functor $\otimes:Spaces\times \C\ra \C $
such that:
\begin{enumerate}
\item there are isomorphisms $\Delta[0]\otimes X\cong X$ natural in
$X$;
\item there are coherent isomorphisms $(K\times
L)\otimes X \cong K\otimes(L\otimes X)$ natural in
 $K,L\in Spaces$ and $X\in \C$ (the explicit coherence diagrams can be found
in~\cite[Definition~4.1.6]{99h:55031});
\item for any $X\in \C$,
$-\otimes X:Spaces\ra \C$ has a right adjoint $\text{map}(X,-)$.
\end{enumerate}

Note that if $\C$ is closed under taking colimits, then condition
(3) can be replaced by the requirement:
\begin{enumerate}
\item[(3')]  for any $X\in \C$, $-\otimes X:Spaces\ra \C$ commutes with arbitrary colimits.
\end{enumerate}

Thus what we call in this paper a left action of the category
$Spaces$ is a special case of what is called a left
$Spaces$-module structure on $\C$
in~\cite[Definition~4.1.6]{99h:55031}, namely one in which the
tensor $-\otimes X$ has a right adjoint. }\end{point}

\begin{point}\label{pt simpfun}{\rm
An enriched functor between categories enriched over $Spaces$ is
called simplicial. A simplicial functor between left actions on
$\C$ and  $\D$ can be described explicitly as a functor
$\phi:\C\ra \D$  together with a natural transformation
$\mu:K\otimes \phi(X)\ra \phi(K\otimes X)$ which satisfies the
following coherence conditions:
\begin{enumerate}
\item the composition $\Delta[0]\otimes
\phi(X)\xra{\mu}\phi( \Delta[0]\otimes X)\cong \phi(X)$ is the
natural isomorphism given in~\ref{pt leftsimpact}(1);
\item for any $X\in \C$ and $K,L\in Spaces$, the following diagram commutes, where the indicated isomorphisms
are induced by the coherent isomorphism in~\ref{pt
leftsimpact}(2):
\[\xymatrix{
 (K\times L)\otimes \phi(X)\rrto^{\mu}\dto^{\cong}& & \phi\big((K\times L)\otimes X\big)\dto^{\cong}\\
K\otimes \big(L\otimes \phi(X)\big)\rto^{K\otimes \mu} & K\otimes
\phi(L\otimes X)\rto^{\mu} & \phi\big(K\otimes (L\otimes X)\big)
}\]
\end{enumerate}
}\end{point}

\begin{point}{\rm
A right action of $Spaces$ on $\C$ is an
 enrichment of $\C$ over
$Spaces$, where the space of morphisms has a right adjoint in the
second variable. A right action of $Spaces$ on $\C$ can be
described explicitly as a functor
 $\text{hom}:Spaces^{op}\times \C\ra \C$ such that:
\begin{enumerate}
\item[(1d)] there are isomorphisms $X\cong\text{hom}(\Delta[0], X)$ natural in $X$;
\item[(2d)] there are coherent isomorphism $\text{hom}(K\times L, X)\cong
\text{hom}\big(K,\text{hom}(L,X)\big)$  natural in $X\in \C$ and
$K,L\in Spaces$;
\item[(3d)] for any $X\in \C$, $\text{hom}(-, X):Spaces^{op}\ra \C$ has a left adjoint
$\text{map}(-,X)$.
\end{enumerate}

Note that $\text{hom}$ is a right action of $Spaces$ on $\C$ if
and only if $\text{hom}^{op}$ is a left action of $Spaces$ on
$\C^{op}$.

A  simplicial functor between right actions of $Spaces$ on $\C$
and $\D$ can be described in a explicit way as a functor
$\phi:\C\ra \D$ together with a natural transformation
$\phi\big(\text{hom}_{\C}(K,Y)\big)\ra
\text{hom}_{\D}\big(K,\phi(X)\big)$ satisfying coherence
conditions analogous to those in~\ref{pt simpfun}. }\end{point}

\begin{point}\label{pt simpstruct}{\rm
A simplicial structure on $\C$ is an enrichment of $\C$ over
$Spaces$ where the space of morphisms has a left adjoint in the
first variable and a right adjoint in the second variable. If such
an enrichment has been fixed, $\C$ is called a {\it simplicial
category}.

A simplicial structure can be described explicitly in two dual
ways; first, as a  left action $\otimes:Spaces\times \C\ra \C$
which satisfies the extra condition:
\begin{enumerate}
\item[(4)] The functor $K\otimes -:\C\ra \C $ has a right adjoint
$\text{hom}(K,-)$ for any $K$.
\end{enumerate}
In this case the assignment $K \mapsto\text{hom}(K,X)$ can be
extended to a functor $Spaces^{op}\ra \C$ which is the right
adjoint to $\text{map}(-,X)$. In particular $\text{hom}$ defines a
right action of $Spaces$ on $\C$.

Dually a simplicial structure on $\C$ can also be  described as a
right action $\text{hom}(-,-):Spaces^{op}\times \C\ra \C$  for
which:
\begin{enumerate}
\item[(4d)] $\text{hom}(K,-):\C\ra \C $ has a left adjoint $K\otimes -$,
for any space $K$.
\end{enumerate}
}\end{point}

\begin{point}\label{pt hosimpstruct}{\rm
If instead of $Spaces$ we take the homotopy category
$\text{Ho}(Spaces)$ in the above explicit descriptions of actions,
then we get the notions of   homotopy (left/right) actions of
$Spaces$ on $\C$ and simplicial functors.  A  weak simplicial
structure on $\C$ is, by definition, a homotopy action. If such an
action is fixed, $\C$ is called a {\it weak simplicial category}.
}\end{point}

\begin{point}\label{pt homsimpfun}{\rm
If $\otimes:\text{Ho}(Spaces)\times \C\ra \C$ is a homotopy left
action, then the set of morphisms $\text{mor}_{\C}(X,Y)$ can be
identified with the set of connected components
$\pi_{0}\text{map}(X,Y)$. It follows that a simplicial functor
between homotopy left actions on $\C$ and $\D$  induces a natural
transformation $\text{map}_{\C}(X,Y)\ra
\text{map}_{\D}\big(\phi(X),\phi(Y)\big)$ which on the set of
connected components coincides with~$\phi$. The same applies to
simplicial functors of homotopy right actions. }\end{point}

\begin{point}\label{pt replaceman}{\rm
Let $(\C,W)$ and $(\D,V)$ be categories with weak equivalences
whose homotopy categories exist. The total left derived functor of
$\Phi:\C\ra \D$ is  the left Kan extension~\cite{MR2001j:18001} of
$\pi_{\C}:\C\ra \text{Ho}(\C)$ applied to the composition
 of $\Phi$ and $\pi_{\D}:\D\ra\text{Ho}(\D)$.
Explicitly, it is a functor $L\Phi:\C\ra \text{Ho}(\D)$ with a
natural transformation $L\Phi\ra \pi_{D}\Phi$ which is terminal
among natural transformations with the range  $\pi_{D}\Phi$ and
whose domains are functors factoring through $\pi_{\C}:\C\ra
\text{Ho}(\C)$.

A {\it left replacement} for $\Phi:\C\ra \D$ is a functor $Q:\C\ra
\C$ and a natural morphism $QX\ra X$ in $\C$ which is a weak
equivalence such that $\pi_{D}\Phi Q\ra \pi_{D} \Phi$  is the
total left derived functor of $\Phi$.

The total right derived functor of $\Phi:\C\ra \D$ is the right
Kan extension~\cite{MR2001j:18001} of $\pi_{\C}:\C\ra \text{\rm
Ho}(\C)$ applied to the composition of
$\pi_{\D}:\D\ra\text{Ho}(\D)$ and $\Phi$. A {\it right
replacement} for $\Phi$ is a functor $R:\C\ra \C$ and a natural
morphism $X\ra  RX$ in $\C$ which is a weak equivalence such that
$\pi_{D}\Phi\ra \pi_{D} \Phi R$  is the total right derived
functor of $\Phi$. }\end{point}

\begin{point}\label{pt localized}{\rm
We say that $(\C,W)$ can be {\it left localized} if, for any small
category $I$, the homotopy category of $\text{Fun}(I,\C)$ exists
and the functor $\text{colim}_{I}:\text{Fun}(I,\C)\ra \C$ has a
left replacement. In particular  the total left derived functor of
$\text{colim}_{I}$ exists and is denoted by
$\text{hocolim}_{I}:\text{Fun}(I,\C)\ra \text{Ho}(\C)$. It is
given by the formula  $\text{\rm
hocolim}_{I}F=\pi_{\C}(\text{colim}_{I} QF)$.

We say that $(\C,W)$ can be {\it right localized} if, for any
small category $I$, the homotopy category of  $\text{Fun}(I,\C)$
exists and the functor $\text{lim}_{I}:\text{Fun}(I,\C)\ra \C$ has
a right replacement. In particular   the total right derived
functor of $\text{lim}_{I}$ exists and is denoted by
$\text{holim}_{I}:\text{Fun}(I,\C)\ra \text{Ho}(\C)$. }\end{point}

\begin{point}{\rm
 Assume that $(\C,W)$ can be left localized. By the universal
property, for any $F:J\ra \C$ and $f:I\ra J$, there is a morphism
$\text{hocolim}_{I}f^{\ast}F\ra \text{hocolim}_{J}F$ in
$\text{Ho}(\C)$  natural with respect to $f$ and $F$. This
morphism can be described in terms of left replacements $Q_{I}$
and $Q_{J}$  of $\text{colim}_{I}$ and $\text{colim}_{J}$.  It is
represented by the following sequence of morphisms in $\C$:
\[
\text{colim}_{I}(Q_{I}f^{\ast}F)\xla{\sim}\text{colim}_{I}(Q_{I}f^{\ast}Q_{J}F)\ra
\text{colim}_{I}(f^{\ast}Q_{J}F)\ra\text{colim}_{J}(Q_{J}F)
\]
In the case $J$ is the category with one object and one morphism,
this leads to a functor:
\[
\otimes_{l}:Cat\times \C\ra \text{Ho}(\C)\ \ \ \ I\otimes_{l}
X:=\text{hocolim}_{I}cX
\]

Dually, when $(\C,W)$ can be right localized, by the universal
property, for any $F:J\ra \C$ and $f:I\ra J$, there is a morphism
$\text{holim}_{J}F\ra \text{holim}_{I}f^{\ast}F$ in
$\text{Ho}(\C)$  natural with respect to $f$ and $F$. In the case
$J$ is the category with one object and one morphism, this leads
to a functor:
\[\text{\rm rhom}:Cat^{op}\times \C\ra \text{\rm Ho}(\C)\ \ \ \
\text{\rm rhom}(I,X):=\text{holim}_{I}cX \] }\end{point}

\begin{point}\label{pt nathocolim}{\rm
Assume that  $(\C,W)$ and $(\D,V)$ are categories  that can be
left localized.  Let $I$ be a small category, $F:I\ra \C$ be a
functor,  and $Q_{\C}$ and $Q_{\D}$ be left replacements of
$\text{colim}_{I}:\text{Fun}(I,\C)\ra\C$ and
$\text{colim}_{I}:\text{Fun}(I,\D)\ra\D$ respectively. If
$\Phi:\C\ra \D$ is a functor that sends weak equivalences in $\C$
to weak equivalences in $\D$, then the following sequence of
morphisms in $\D$:
\[\text{colim}_{I}Q_{\D}\Phi F\xla{\sim}\text{colim}_{I}Q_{\D}\Phi Q_{\C}F\ra
\text{colim}_{I}\Phi Q_{\C}F\ra \Phi(\text{colim}_{I}Q_{\C}F)\]
leads to a morphism $\text{hocolim}_{I}\Phi F\ra
\Phi(\text{hocolim}_{I}F)$ in $\text{Ho}(\D)$. This morphism is
natural in $F$ and does not depend on the choice of the left
replacements $Q_{\C}$ and $Q_{\D}$.

Dually, let $(\C,W)$ and $(\D,V)$ be categories that can be right
localized. Assume that $\Phi:\C\ra \D$  sends weak equivalences in
$\C$ to weak equivalences in $\D$. As in the case of hocolims, for
any small category $I$ and any functor $F:I\ra \C$, there is a
natural  morphism $\Phi(\text{holim}_{I}F)\ra \text{holim}_{I}\Phi
F$ in $\text{Ho}(\D)$. }\end{point}

\section{Simplex Categories}
\label{sec simplexcat}
To understand  actions  of $Spaces$ one can try to look at the
problem locally and study functors indexed by  simplex categories
(\ref{pt simplexcat}). Such functors are called representations of
the indexing space. What is so special about simplex categories
that enables us to say something about their representations?

Simplex categories have one important advantage over arbitrary
small categories of which we give now an elementary illustration.
 Let $g:I\ra X$ be a functor of
small categories and $G:I\ra \C$ be a functor. Let us look for
conditions under which the operation $G\mapsto g^{k}G$ commutes
with the base change. Explicitly, consider a pull-back of small
categories:
\[
\xymatrix{
P\rto^{ f_{1} }\dto_{g_{1}} & I\dto ^{g}\\
J\rto^{f} & X }
\]
For any $j\in J$, there is a functor of over-categories
$g_{1}\downarrow j\ra g\downarrow f(j)$. These functors induce a
natural transformation:
\[
(g_{1}^{k}f_{1}^{\ast}G)(j)=\text{colim}_{g_{1}\downarrow j} \,
f_{1}^{\ast}G \ra \text{colim}_{g\downarrow
f(j)}G=(f^{\ast}g^{k}G)(j)
\]
We ask: when is this natural transformation an isomorphism? If
$G=g^{\ast}F$, for some  $F:X\ra \C$, this natural  transformation
is of the form:
\[
{g_{1}}^{k}{g_{1}}^{\ast}f^{\ast}F=
{g_{1}}^{k}{f_{1}}^{\ast}g^{\ast}F\ra f^{\ast} g^{k}g^{\ast}F
\]
We could then apply $f^{k}$ to it and ask: when is
$(fg_{1})^{k}(fg_{1})^{\ast}F\ra f^{k}f^{\ast}g^{k}g^{\ast}F$ an
isomorphism? In particular, when is the operation $F\mapsto
f^{k}f^{\ast}g^{k}g^{\ast}F$ symmetric with respect to $f$ and
$g$?

Let us look at two examples:
\begin{example}{\rm
First a negative case. Consider a pull-back  of categories:
\[
\xymatrix{
\emptyset\rto^{f_{1}}\dto_{g_{1}} & 1\dto^{g}\\
0\rto^(0.3){f} & \{0\ra 1\}\hspace{7mm} }
\]
and a functor $\{A\ra B\}:\{0\ra 1\}\ra Spaces$. It can be
verified that:
\[
g^{k}g^{\ast}\{A\ra B\}=\{\emptyset\ra B\}\ \ \hbox{\rm and}\ \ \
f^{k}f^{\ast}\{A\ra B\}=\{A\stackrel{id}{\ra}A\}
\]
In particular $f^{k}f^{\ast}g^{k}g^{\ast}\{A\ra B\}=
\{\emptyset\ra\emptyset\}$ does not coincide with the functor
$g^{k}g^{\ast}f^{k}f^{\ast}\{A\ra B\}=\{\emptyset\ra A\}$ in
general.

\medskip

 Here is a positive case. Consider a  pull-back of
categories:
\[
\xymatrix{
\emptyset\rto^{f_{1}}\dto_{g_{1}}  & 1\dto^{g}\\
0\rto^(0.2){f}  & \{0\ra 01\la 1\}\hspace{16mm} }\] and a functor
$\{A\ra B\la C\}:\{0\ra 01\la 1\}\ra Spaces$. In this case:
\[
g^{k}g^{\ast}\{A\ra B\la C\}= \{\emptyset\ra C\stackrel{id}{\la}
C\}\]
\[f^{k}f^{\ast}\{A\ra B\la C\}= \{A\stackrel{id}{\ra} A\la \emptyset\}
\]
In particular:
\[f^{k}f^{\ast}g^{k}g^{\ast}\{A\ra B\la C\}= \{\emptyset\ra \emptyset\la \emptyset\}\]
\[g^{k}g^{\ast}f^{k}f^{\ast}\{A\ra B\la C\}= \{\emptyset\ra \emptyset\la \emptyset\}\]
which agrees with $(gf_{1})^{k}(gf_{1})^{\ast}\{A\ra B\la C\}$.
}\end{example}

What is so special about the second example? It turns out that it
is an example of functors indexed by simplex categories. The
reason this does not work in general is the fact that arbitrary
categories can have vastly different over-categories. Simplex
categories on the other hand look the same locally. For any
simplex $\sigma:\Delta[n]\ra K$, the functor from the
over-category $\mathbf K\downarrow \sigma\ra \mathbf K$ is
isomorphic to the functor induced by the map $\sigma: \Delta[n]
\ra  K$. Thus, locally, all simplex categories look like standard
simplices.

\begin{proposition}\label{prop pullback}
Let the following be a pull-back square of spaces:
\[\xymatrix{
P\rto^{f_{1}}\dto_{g_{1}}  & A\dto^{g}\\
B\rto^{f} & D }
\] For any functor $G:\mathbf A\ra \C$, the natural
transformation ${\mathbf g_{1}}^{k}{\mathbf f_{1}}^{\ast}G\ra
{\mathbf f}^{\ast}{\mathbf g}^{k}G$ is an isomorphism. In
particular, for any $F:\mathbf D\ra \C$, the natural
transformations ${\mathbf g_{1}}^{k}{\mathbf g_{1}}^{\ast}{\mathbf
f}^{\ast}F\ra {\mathbf f}^{\ast}{\mathbf g}^{k}{\mathbf
g}^{\ast}F$ and $({\mathbf f}{\mathbf g_{1}})^{k}({\mathbf
f}{\mathbf g_{1}})^{\ast}F\ra {\mathbf f}^{k}{\mathbf
f}^{\ast}{\mathbf g}^{k}{\mathbf g}^{\ast}F$ are isomorphisms and
the operation $F\mapsto {\mathbf f}^{k}{\mathbf f}^{\ast}{\mathbf
g}^{k}{\mathbf g}^{\ast}F$ is symmetric with respect to $f$ and
$g$.
\end{proposition}

\begin{proof}
Note that for any simplex $\sigma:\Delta[n]\ra B$, the
over-categories ${\mathbf g_{1}}\downarrow\sigma$ and ${\mathbf
g}\downarrow f(\sigma)$ can be identified with the simplex
category of the pull-back:
\[
\text{lim}(\Delta[n]\stackrel{\sigma}{\ra}B\stackrel{g_{1}}{\la}P)\]
Moreover the functor ${\mathbf g_{1}}\downarrow\sigma\ra {\mathbf
g}\downarrow f(\sigma)$ is an isomorphism.
\end{proof}

\section{Bounded functors}\label{sec bounded}
Simplex categories have another advantage. There are a lot of
geometric constructions that can be performed on spaces. These
constructions often translate well to operations on
representations. For example, fix a functor $H:I\ra Spaces$. The
simplex category of $\colim{I}{H}$ can be identified with the
colimit of $\mathbf H$ in  the category of small categories. A
covariant functor $F$ indexed by  $\colim{I}{\mathbf H}$ consists
of the following data:
\begin{enumerate}
\item for every object  $i$ in $I$, a functor $F_{i}:\mathbf H(i)\ra \C$;
\item for every morphism $\alpha:i\ra j$ in $I$, $F_{i}:\mathbf H(i)\ra
\C$ should coincide with the composition $\mathbf
H(i)\stackrel{\mathbf H(\alpha)}{\lra}\mathbf
H(j)\stackrel{F_{j}}{\lra} \C$.
\end{enumerate}
This decomposition yields an isomorphism
$\text{colim}_{\text{colim}_{I}\mathbf H}{F}\cong
\colim{I}{\colim{\mathbf H(i)}{F_{i}}}$, see
\cite[Proposition~8.2]{MR2002k:55026}.

By applying the above procedure inductively, we can calculate
colimits using the cell decomposition of the indexing space.
Diagrams indexed by the standard simplices $\mathbf \Delta [n]$
play the role of fundamental building blocks in this process.
Unfortunately representations of $\Delta[n]$ are complicated. For
example, the simplex category of  $\Delta[0]$ is equivalent to the
category of finite ordered sets and so its representations are
cosimplicial and simplicial objects, depending on variance. For
that reason we need to make additional assumptions on the functors
considered. We restrict our attention to diagrams that reflect the
geometry of the indexing space more directly. A functor $F$,
indexed by a simplex category $\mathbf A$, is called {\em bounded}
\cite[Definition 10.1]{MR2002k:55026} if it assigns an isomorphism
to any morphism of the form $s_{i}:\Delta[n+1]\ra \Delta[n]$ in
$\mathbf A$.

If $S$ denotes the set of all these degeneracy morphisms in
$\mathbf A$, then a bounded functor is nothing else but a functor
indexed by the localized category $\mathbf A[S^{-1}]$. Depending
on the variance, the categories of bounded diagrams are denoted by
$\text{Fun}^{b}(\mathbf A,\C)$ and $\text{Fun}^{b}(\mathbf
A^{op},\C)$.

An important observation is that the degeneracy morphisms in
$\mathbf A$ do not form any essential loop. Any bounded functor is
naturally isomorphic to a functor that assigns identities to
morphisms in $S$ \cite[Proposition 10.3]{MR2002k:55026}. This can
be therefore assumed about any considered bounded functor.

\begin{example}\label{ex nonhoinv}{\rm
The ``triviality'' of $S$ does not imply in general that the
localization functor
 $\mathbf A\ra \mathbf A[S^{-1}]$ is a weak equivalence (see~\ref{pt wecat}). Let:
\[
A=\text{colim}(\Delta[3]/\Delta[2]\xleftarrow{\pi}
\Delta[3,1]\mono \Delta[3])
\]
where $\pi:\Delta[3,1]\ra \Delta[3]/\Delta[2]$ is the composition
of  $\Delta[3,1]\mono \Delta[3]$ and the quotient map
$\Delta[2]\xra{d_{1}} \Delta[3]\ra \Delta[3]/\Delta[2]$. Although
$A$ is contractible, the category $\mathbf A[S^{-1}]$ is not. The
nerve $N(\mathbf A[S^{-1}])$ is weakly equivalent to $S^{2}$.

Let $A=\Delta[n]/\partial\Delta[n]$. For $n>1$, the category
$\mathbf A[S^{-1}]$ is equivalent to $\ast\ra \ast$ (the category
with two objects and one non-identity morphism). }\end{example}

It is clear that, for any map of spaces $f:A\ra B$, ${\mathbf
f}^{\ast}:\text{Fun}({\mathbf B},\C)\ra \text{Fun}({\mathbf
A},\C)$ takes bounded functors to bounded functors. The same is
also true for the left Kan extension ${\mathbf
f}^{k}:\text{Fun}({\mathbf A},\C)\ra \text{Fun}({\mathbf B},\C)$
\cite[Theorem 33.1]{MR2002k:55026}. It follows that the functors
${ \mathbf f}^{k}:\text{Fun}^{b}({\mathbf A},\C)\rightleftarrows
\text{Fun}^{b}({\mathbf B},\C) :{\mathbf f}^{\ast}$ form an
adjoint pair.

\section{Subdivision}\label{sec subdivision}
In this section we study  bounded functors indexed by barycentric
subdivisions. Recall that the subdivision of $A$ is given by the
nerve $N(\mathbf A)$ of the simplex category of a space~$A$
(see~\ref{pt subdiv}).

A map of spaces $f:A\ra B$ is called {\em reduced} if it sends
non-degenerate simplices in $A$ to non-degenerate simplices in $B$
\cite[Definition 12.9]{MR2002k:55026}. Reduced maps can be
characterized using a lifting test: the map $f:A\ra B$ is reduced
if and only if in any commutative diagram of the form:
\[\xymatrix{
\Delta[n+1]\dto_{s_{i}}\rto & A\dto^{f}\\
\Delta[n]\rto & B }
\]
there is a lift, i.e. a morphism $\Delta[n]\ra A$ for which the
resulting diagram with five arrows commutes~\cite[Proposition
18.5]{MR2002k:55026}. It follows that a pull-back of a reduced map
is also reduced.

The subdivision operation has two crucial properties. First:
\begin{proposition}[{\cite[Example 12.10]{MR2002k:55026}}]\label{prop simpl}
For any map of spaces $f:A\ra B$, the induced map $N(\mathbf
f):N(\mathbf A)\ra N(\mathbf B)$ is reduced. \hfill{\qed}
\end{proposition}

Second, the subdivision $N(\mathbf A)$ contains enough
combinatorial information to induce a non-trivial simplicial
structure on $\text{Fun}^{b}\big(\mathbf N(\mathbf A),\C\big)$
(the non-triviality will follow from Corollary~\ref{col
kanhomclas}). Here is the description of this structure. Choose a
space $K$ and denote by $\pi:N(\mathbf{K}\boxtimes \mathbf{A})\ra
N(\mathbf A)$ the subdivision of the projection
  $p_{K}:K\times
A\ra A$  onto the second factor. Define:
\[
K\otimes F:= {\mathbf \pi}^{k}{\mathbf \pi}^{\ast}F
\]
This construction is natural in $F$ and $K$ and hence defines a
functor:
\[
\otimes :Spaces\times \text{Fun}^{b}\big(\mathbf N(\mathbf
A),\C\big) \ra \text{Fun}^{b}\big(\mathbf N(\mathbf A),\C\big)
\]

\begin{example}{\rm
In general, we do not know any explicit formula describing
$K\otimes F$ in terms of the values of $F$. However if $F$ is a
composition of $G:\mathbf A^{op}\ra \C$ and $\epsilon: {\mathbf
N(\mathbf A)}\ra \mathbf A^{op}$, $\epsilon(\sigma_{n}\ra
\sigma_{n-1}\ra\cdots\ra\sigma_{0}):= \sigma_{n}$, then $K\otimes
F$ is given by:
\[
(\sigma_{n}\ra \sigma_{n-1}\ra\cdots\ra\sigma_{0})\mapsto
\coprod_{K_{|\sigma_{n}|}}G(\sigma_{n})
\]
This formula should be compared with the simplicial structure
introduced by Quillen on simplicial objects in \cite[II.1
Proposition~2]{MR36:6480}.}\end{example}

\begin{proposition}\label{prop simplicial}
For any space $A$, the functor $\otimes$ defines a simplicial
structure on $\text{Fun}^{b}\big({\mathbf N(\mathbf A)},\C\big)$
(see~\ref{pt simpstruct} and \ref{pt leftsimpact}).
\end{proposition}

\begin{proof}
We need to check four conditions (see~\ref{pt simpstruct}). It is
clear that the natural morphism $\Delta[0]\otimes
F=\text{id}^{k}\text{id}^{\ast}F\ra F$ is an isomorphism and so
the first requirement is satisfied.

Consider two spaces $K$ and $L$ and a pull-back square:
\[
\xymatrix{
N({\mathbf K}\boxtimes {\mathbf L}\boxtimes {\mathbf A})\rto^{\bar\pi_{K}}\dto_{\bar \pi_{L}} & N({\mathbf L}\boxtimes {\mathbf A})\dto^{\pi_{L}}\\
N({\mathbf K}\boxtimes {\mathbf A})\rto^{\pi_{K}} & N({\mathbf A})
}
\]
According to Proposition~\ref{prop pullback}, for any $F:\mathbf
N(\mathbf A)\ra \C$, the transformation:
\[
 (K\times L)\otimes F=(\pi_{K}{\bar \pi_{L}})^{k} (\pi_{K}{\bar
\pi_{L}})^{\ast}F\ra
(\pi_{K})^{k}(\pi_{K})^{\ast}(\pi_{L})^{k}(\pi_{L})^{\ast}F =
K\otimes (L\otimes F)
\]
is an isomorphism. Coherence of these isomorphisms follows
directly from the coherence of products in $Spaces$, which proves
the second requirement.

To prove  condition (3'), we need to show that, for a fixed $F\in
\text{Fun}^{b}\big(\mathbf N(\mathbf A),\C\big)$,  the functor
$-\otimes F$ commutes with arbitrary colimits. Let $H:I\ra Spaces$
be a functor. Observe that the map $\text{colim}_{I}(H\times A)\ra
(\text{colim}_{I}H)\times A$ is an isomorphism of spaces and so is
the map $\text{colim}_{I}N({\mathbf H}\boxtimes {\mathbf A})\ra
N\big(( \text{colim}_{I}\mathbf H)\boxtimes{\mathbf A}\big)$ as
the subdivision commutes with arbitrary colimits (see~\ref{pt
subdiv}). Next, choose a simplex $\sigma:\Delta[n]\ra N(\mathbf
A)$ and define:
\[
P(i):=\text{lim}\Big(\Delta[n]\stackrel{\sigma}{\ra} N({\mathbf
A}) \stackrel{\pi_{i}}{\la} N\big({\mathbf H(i)}\boxtimes {\mathbf
A}\big)\Big)
\]
where $\pi_{i}$ is the subdivision of the projection
$pr_{2}:H(i)\times A\ra A$. The colimit of $P$ can be identified
with the pull-back:
\[\text{colim}_{I}P\cong\text{lim}\Big(\Delta[n]\stackrel{\sigma}{\ra}
 N(\mathbf A)\stackrel{\pi}{\la} N\big((\text{colim}_{I}\mathbf H)\boxtimes {\mathbf A}\big)\Big)
\]
These isomorphisms yield the following identifications:

\begin{align*}
\big((\text{colim}_{I}H)\otimes F\big)(\sigma) & =
(\pi^{k}\pi^{\ast}F)(\sigma) \cong
\text{colim}_{\text{colim}_{I}\mathbf P}F \cong \text{colim}_{i\in
I}\text{colim}_{\mathbf P(i)}F \\
& =\text{colim}_{i\in
I}\big((\pi_{i})^{k}(\pi_{i})^{\ast}F\big)(\sigma)=
\text{colim}_{I}(H\otimes F)(\sigma)
\end{align*}

Note that for a fixed $K$, the functor $K\otimes -$ has a right
adjoint as it is a composition of two functors that admit right
adjoints (here we use our general assumption that $\C$ has all
limits, see~\ref{pt limcolim}). This proves the fourth
requirement.
\end{proof}

Taking the subdivision in Proposition~\ref{prop simplicial} is
crucial. In a similar way one could define a simplicial structure
on $\text{Fun}^{b}(\mathbf A,\C)$ by $K\otimes F:=({\mathbf
p_{K}})^{k}({\mathbf p_{K}})^{\ast}F$, where $p_{K}:K\times A\ra
A$ is the projection onto the second factor. This however does not
give anything interesting as the obtained structure is trivial.

\begin{point}{\bf Naturality}.
\label{naturality}{\rm
How functorial is this simplicial structure on
$\text{Fun}^{b}({\mathbf N(\mathbf A)},\C)$? Consider the
following pull-back of spaces associated with a map $f:A\ra B$ and
a space~$K$:
\[\xymatrix{
N(\mathbf K\boxtimes \mathbf A)\dto_{N( \mathbf{id}\boxtimes \mathbf{f})} \rto^(0.58){\pi} & N(\mathbf A)\dto^{N(\mathbf f)} \\
N(\mathbf K\boxtimes \mathbf B)\rto ^(0.58){\pi_{1}}  & N(\mathbf
B) }\] where the horizontal maps are the subdivisions of the
projections. For any bounded diagram $F\in\text{Fun}^{b}(\mathbf
N(\mathbf B),\C)$, there is a natural transformation:
\[
K\otimes \big(\mathbf N(\mathbf
f)^{\ast}F\big)=\pi^{k}\pi^{\ast}\mathbf N(\mathbf f)^{\ast}F\ra
\mathbf N(\mathbf f)^{\ast} (\pi_{1})^{k}(\pi_{1})^{\ast}F=\mathbf
N(\mathbf f)^{\ast}(K\otimes F)
\]
which is an isomorphism by Proposition~\ref{prop pullback}. Note
that this isomorphism is natural with respect to the variables $K$
and $F$. Together with this isomorphism the functor $\mathbf
N(\mathbf f)^{\ast}:\text{Fun}^{b}(\mathbf N(\mathbf B),\C)\ra
\text{Fun}^{b}(\mathbf N(\mathbf A),\C)$ becomes a simplicial
functor (see~\ref{pt simpfun}).

For any $G\in\text{Fun}^{b}(\mathbf N(\mathbf A),\C)$, there is a
natural transformation:
\[
\mathbf N(\mathbf{id}\boxtimes
\mathbf{f})^{k}\pi^{\ast}G\ra(\pi_{1})^{\ast}\mathbf N(\mathbf
f)^{k}G
\]
which is an isomorphism again by Proposition~\ref{prop pullback}.
We can apply $(\pi_{1})^{k}$ to it and get a natural isomorphism:
\[\mathbf N(\mathbf f)^{k}(K\otimes G)=
(\pi_{1})^{k}\mathbf N(\mathbf{id}\boxtimes
\mathbf{f})^{k}\pi^{\ast}G\ra (\pi_{1})^{k}(\pi_{1})^{\ast}\mathbf
N(\mathbf f)^{k}G= K\otimes \big(\mathbf N(\mathbf f)^{k}G\big)
\]
where we identified $K\otimes G=\pi^{k}\pi^{\ast}G$ and used the
commutativity of the above pull-back square. It is not difficult
to see that this isomorphism is natural with respect to $K$
and~$G$. The functor $\mathbf N(\mathbf
f)^{k}:\text{Fun}^{b}\big(\mathbf N(\mathbf A),\C\big)\ra
\text{Fun}^{b}\big(\mathbf N(\mathbf B),\C\big)$, together with
the inverse of this isomorphism, becomes a simplicial functor.

Let $\phi:\C\ra \D$ be a functor and $K$ be a simplicial set.
Consider the following pull-back of spaces:
\[\xymatrix{
P_{\sigma}\rto\dto & N({\mathbf K}\boxtimes {\mathbf A})\dto^{\pi}\\
\Delta[n]\rto^{\sigma} & N(\mathbf A) }
\]
By the universal property of the colimit construction, for any
$F:\mathbf N(\mathbf A)\ra \C$, there is a morphism in $\D$:
\[
\big(K\otimes (\phi F)\big)(\sigma)=\text{colim}_{\mathbf
P_{\sigma}}(\phi F \pi)\ra \phi(\text{colim}_{\mathbf
P_{\sigma}}F\pi)=\phi(K\otimes F)(\sigma)
\]
These morphisms form a natural transformation $K\otimes (\phi
F)\ra \phi(K\otimes F)$ functorial with respect to the variables
$K$ and $F$. With this natural transformation, the functor
$\phi_{\ast}: \text{Fun}^{b}\big(\mathbf N(\mathbf A),\C\big)\ra
\text{Fun}^{b}\big(\mathbf N(\mathbf A),\D\big)$ becomes a
simplicial functor.}
\end{point}

\section{Homotopy theory of bounded functors}
\label{sec hobounded}
Simplex categories are also rather special from the homotopy point
of view.

\begin{theorem}[{\cite[Theorem 21.1]{MR2002k:55026}}]\label{thm modelmain}
Let $\M$ be a model category. The following describes a model
structure on $\text{\rm Fun}^{b}(\mathbf A,\M)$:
\begin{itemize}
\item $\phi:F\ra G$ is a weak equivalence (fibration) if,
for any simplex $\sigma\in A$, $\phi_{\sigma}:F(\sigma)\ra
G(\sigma)$ is a weak equivalence (fibration) in $\M$;
\item $\phi: F\ra G$ is an (acyclic) cofibration if, for any non-degenerate
simplex $\sigma:\Delta[n]\ra A$, the morphism:
\[
\text{\rm colim}\big(F(\sigma)\la \text{\rm
colim}_{\mathbold{\partial \Delta}[n]}F\ra \text{\rm
colim}_{\mathbold{\partial \Delta}[n]}G\big)\ra G(\sigma)
\]
is an (acyclic) cofibration in $\M$.\hfill{\qed}
\end{itemize}
\end{theorem}

The description of (acyclic) cofibrations in Theorem~\ref{thm
modelmain} can be used by induction on the cell attachment to
show:

\begin{proposition}[{\cite[Proposition 20.1]{MR2002k:55026}}]\label{prop ext}
If  $\phi:F\ra G$ is an (acyclic) cofibration in $\text{\rm
Fun}^{b}(\mathbf A,\M)$, then, for any subspace $X\subset A$, the
morphism:
\[
\text{\rm colim}\big(\text{\rm colim}_{\mathbf X}
G\stackrel{\phi}{\la}\text{\rm colim}_{\mathbf X}F\ra \text{\rm
colim}_{\mathbf A}F\big)\hookrightarrow \text{\rm colim}_{\mathbf
A}G
\]
is an (acyclic) cofibration in $\M$.\hfill{\qed}
\end{proposition}

How functorial is this model structure?

\begin{proposition}\label{prop basicprophocolim}
Let $f:A\ra B$ be a map of $Spaces$.
\begin{enumerate}
\item  If $\phi$ is an (acyclic) cofibration  in $\text{\rm
Fun}^{b}(\mathbf A,\M)$, then ${\mathbf f}^{k}\phi$ is an
(acyclic) cofibration  in $\text{\rm Fun}^{b}(\mathbf B,\M)$
{\cite[Theorem 11.2.]{MR2002k:55026}}. In particular $\text{\rm
colim}_{\mathbf A}\phi$ is an (acyclic) cofibration in $\M$. \item
If $f$ is reduced and $\phi$ is an (acyclic) cofibration in
$\text{\rm Fun}^{b}(\mathbf B,\M)$, then so is ${\mathbf
f}^{\ast}\phi$ in $\text{\rm Fun}^{b}(\mathbf A,\M)$.
\item If $\phi$ is
an (acyclic) cofibration in $\text{\rm Fun}^{b}(\mathbf N(\mathbf
B),\M)$, then so is $\mathbf N(\mathbf f)^{\ast}\phi$ in the
category $\text{\rm Fun}^{b}(\mathbf N(\mathbf
A),\M)$.\hfill{\qed}
\end{enumerate}
\end{proposition}

The model structure given in  Theorem~\ref{thm modelmain} was
introduced to study homotopy invariance of the colimit operation.
The motivation was to explain:
\begin{theorem}[{\cite[Theorem 20.3, Corollary 21.4]{MR2002k:55026}}]
\label{thm hoinvcolim}
Let $\phi:F\ra G$ be a natural transformation in $\text{\rm
Fun}^{b}(\mathbf A,\M)$. If $\phi$ is a weak equivalence between
cofibrant objects, then so is $\text{\rm colim}_{\mathbf A}\phi$.
\hfill{\qed}
\end{theorem}

Theorem~\ref{thm hoinvcolim} implies that the functor
$\text{colim}_{\mathbf A}:\text{Fun}^{b}(\mathbf A,\M)\ra \M$ has
a left replacement (see~\ref{pt replaceman}) given by a cofibrant
replacement $Q$ in $\text{Fun}^{b}(\mathbf A,\M)$. In particular
its total left derived functor denoted by $\text{ocolim}_{\mathbf
A}:\text{Fun}^{b}(\mathbf A,\M)\ra \text{Ho}(\M)$ assigns to a
bounded functor $F$ the object $\text{colim}_{\mathbf A}(QF)$ in
$\M$ \cite[Proposition 14.2]{MR2002k:55026}. The functor
$\text{ocolim}_{\mathbf A}$ coincides with
$\text{hocolim}_{\mathbf A[S^{-1}]}$ (the total left derived
functor of $\text{colim}_{\mathbf A[S^{-1}]}:\text{Fun}(\mathbf
A[S^{-1}],\M)\ra \M$). In general $\text{ocolim}_{\mathbf A}$ is
{\em not} the restriction of  the total left derived functor of
$\text{colim}_{\mathbf A}:\text{Fun}(\mathbf A,\M)\ra \M$. To
construct $\text{hocolim}_{\mathbf A}$  an additional step is
necessary: $A$ needs to be subdivided first. This works in general
for diagrams indexed by any small category $I$ and with values in
any category that admits a left model approximation.

\begin{theorem}[{\cite[Theorem 11.3]{MR2002k:55026}}]\label{thm hocolimconst}
Let $I$ be a small category, $\epsilon:\mathbf N(I)\ra I$ the
forgetful functor $(i_{n}\ra\cdots\ra i_{0})\mapsto i_{0}$ and
$l:\M\leftrightarrows \C:r$ a left model approximation. Then the
following functors form a left model approximation:
\[
l\epsilon^{k}:\text{\rm Fun}^{b}(\mathbf
N(I),\M\big)\rightleftarrows \text{\rm Fun}(I,\C):r\epsilon^{\ast}
\]
Moreover, $F\mapsto l\epsilon^{k}Q(r\epsilon^{\ast}F)$ is a left
replacement for $\text{\rm colim}_{I}:\text{\rm Fun}(I,\C)\ra \C$.
\hfill{\qed}
\end{theorem}

\begin{corollary}\label{col leftlocalization}
If  $(\C,W)$ admits a left model approximation, then it can be
left localized (see~\ref{pt localized}). In particular there is a
functor $\otimes_{l}:Cat\times \C\ra \text{\rm Ho}(\C)$ defined by
the assignment $(I,X)\mapsto \text{\rm hocolim}_{I}cX$.
\end{corollary}

\begin{proof}
According to Theorem~\ref{thm hocolimconst}, any functor $F:I\ra
\C$ can be modified functorially to
$Q_{I}F:=l\epsilon^{k}Q(r\epsilon^*F)$ so that
$\text{colim}_{I}Q_{I}F$ represents $\text{hocolim}_{I}F$ in
$\text{\rm Ho}(\C)$.
\end{proof}

For categories admitting left model approximations, homotopy
colimits can be characterized in an alternative way:

\begin{corollary}\label{col althocolim}
Assume that $(\C,W)$ admits a left model approximation. The
functor $\text{\rm hocolim}_{I}:\text{\rm Ho}\big(\text{\rm
Fun}(I,\C)\big)\ra \text{\rm Ho}(\C)$ is then left adjoint to the
constant diagram functor $c:\text{\rm Ho}(\C)\ra \text{\rm
Ho}\big(\text{\rm Fun}(I,\C)\big)$.
\end{corollary}

\begin{proof}
If $\M$ is   a model category, then $\text{colim}_{\mathbf
K}:\text{Fun}^{b}(\mathbf K,\M)\ra \M$ preserves (acyclic)
cofibrations (Proposition~\ref{prop basicprophocolim}). Thus
by~\cite[Theorem 9.7]{MR96h:55014}, the total derived functors of
the  adjoint functors $\text{colim}_{\mathbf
K}:\text{Fun}^{b}(\mathbf K,\M)\rightleftarrows \M:c$ also form an
adjoint pair $\text{ocolim}_{\mathbf K}:\text{Ho}(
\text{Fun}^{b}\big(\mathbf K,\M)\big)\rightleftarrows
\text{Ho}(\M):c$.

Choose a left model approximation $l:\M\leftrightarrows \C:r$. Let
$F:I\ra \C$ be a functor  and $Y$ an object in $\C$. We have the
following bijections between sets of morphisms respectively in
$\text{Ho}\big(\text{Fun}(I,\C)\big)$, $\text{Ho}\big(\text{\rm
Fun}^{b}(\mathbf N(I),\M)\big)$, $\text{Ho}(\M)$, and
$\text{Ho}(\C)$:
\begin{align*}
[F,cY] &
\stackrel{(a)}{\cong}[r\epsilon^{\ast}F,r\epsilon^{\ast}cY]=
[r\epsilon^{\ast}F,crY]\stackrel{(b)}{\cong}
[\text{ocolim}_{\mathbf N(I)}r\epsilon^{\ast}F, rY] \\
& \stackrel{(c)}{\cong}[L(l)\text{ocolim}_{\mathbf
N(I)}r\epsilon^{\ast}F, Y]= [\text{hocolim}_{I}F,Y]
\end{align*}
where (a) follows from the fact that $l\epsilon^{k}:\text{\rm
Fun}^{b}\big(\mathbf N(I),\M\big)\rightleftarrows \text{\rm
Fun}(I,\C):r\epsilon^{\ast}$ is a left model approximation
(Theorem~\ref{thm hocolimconst}) and \cite[Proposition
5.5]{MR2002k:55026}, (b) is a consequence of adjointness of
$\text{ocolim}_{\mathbf N(I)}$ and $c$, and the last bijection (c)
follows from the following Lemma~\ref{lem adjmodapprox}.
\end{proof}

\begin{lemma}\label{lem adjmodapprox}
If $l:\M\rightleftarrows \C:r$ is a left model approximation, then
the total derived functors induce an adjoint pair $L(l):\text{\rm
Ho}(\M)\rightleftarrows \text{\rm Ho}(\C):r$.
\end{lemma}

\begin{proof}
The lemma is a consequence of two facts. First $r:\text{\rm
Ho}(\C)\ra \text{\rm Ho}(\M)$ is a fully-faithful embedding
\cite[Proposition 5.5]{MR2002k:55026}. Second $L(l)(rX)\ra X$ is
an isomorphism in $\text{\rm Ho}(\C)$, for any $X$.
\end{proof}

\section{Constant functors}\label{sec constatfun}
Assume that $(\C,W)$ admits a left model approximation. To
understand how  this structure leads to an action of $Spaces$ on
$\C$ we are going to study properties of homotopy colimits of
constant functors. The fundamental observation is:
\begin{proposition}\label{prop hoinvop}
If $f:I\ra J$ is a weak equivalence of small categories, then, for
any $X\in \C$, $f\otimes_{l}X:I\otimes_{l}X\ra J\otimes_{l}X$ is
an isomorphism in $\text{\rm Ho}(\C)$.
\end{proposition}

\begin{corollary}\label{col hoinvcf}
If $f:K\ra L$ is a weak equivalence of spaces, then, for any $X\in
\C$, ${\mathbf f}\otimes_{l}X:{\mathbf K}\otimes_{l}X\ra {\mathbf
L}\otimes_{l}X$ is an isomorphism in $\text{\rm Ho}(\C)$.
\hfill{\qed}
\end{corollary}

Corollary~\ref{col hoinvcf} implies that the functor $(K,X)\mapsto
\text{\rm hocolim}_{\mathbf K}cX={\mathbf K}\otimes_{l} X$  is a
composition of  the localization $Spaces\times \C\ra
 \text{\rm Ho}(Spaces)\times\text{\rm
Ho}(\C)$ and a functor which we denote by the same symbol:
\[
\otimes_{l}:\text{\rm Ho}(Spaces)\times \text{\rm Ho}(\C)\ra
\text{\rm Ho}(\C)
\]
The aim of this paper is to understand to what extent
$\otimes_{l}$ is a homotopy action. We start by proving that it
satisfies the first two requirements of 2.16 and will show in
Proposition~\ref{prop mapasadjoint} that property (3) holds as
well.

\begin{corollary}\label{lem prophococons}\hspace{2mm}
\begin{enumerate}
\item There are isomorphisms ${\mathbf \Delta[0]}\otimes_{l}X\cong X$ in
$\text{\rm Ho}(\C)$  natural in $X\in \text{\rm Ho}(\C)$.
\item There are coherent isomorphisms
$({\mathbf K}\boxtimes {\mathbf L})\otimes_{l}X\cong {\mathbf
K}\otimes_{l}({\mathbf L}\otimes_{l} X)$ in $\text{\rm Ho}(\C)$
natural in $X\in \text{\rm Ho}(\C)$ and $K,L \in \text{\rm
Ho}(Spaces)$.
\end{enumerate}
\end{corollary}
\begin{proof}
Part (1) follows from Proposition~\ref{prop hoinvop} and the fact
that the functor between the simplex category  of $\Delta[0]$ and
the category with only one object and one morphism is a weak
equivalence of categories.

Recall that ${\mathbf K}\boxtimes {\mathbf L}\mono {\mathbf
K}\times {\mathbf L}$ is a weak equivalence of categories. Thus by
Proposition~\ref{prop hoinvop}, $({\mathbf K}\boxtimes {\mathbf
L})\otimes_{l}X\ra ({\mathbf K}\times {\mathbf L})\otimes_{l}X$ is
an isomorphism in $\text{Ho}(\C)$. According to~\cite[Theorem
24.9]{MR2002k:55026} there is a further isomorphism in
$\text{Ho}(\C)$:
\[
(\mathbf K\times \mathbf L)\otimes_{l}X=\text{hocolim}_{{\mathbf
K}\times {\mathbf L}}cX \cong \text{hocolim}_{\mathbf
K}\text{hocolim}_{\mathbf L}cX= {\mathbf K}\otimes_{l}({\mathbf
L}\otimes_{l}X)
\]
The coherence of these isomorphisms is straightforward. This
concludes the proof.
\end{proof}

\noindent {\it Proof of~\ref{prop hoinvop}.} Since $\epsilon:
\mathbf N(I)^{op}\ra I,\ (i_{n}\ra\cdots\ra i_{0})\mapsto i_{n}$
is terminal~\cite[Theorem 29.1]{MR2002k:55026}, the morphism
$\text{hocolim}_{\mathbf N(I)^{op}}cX\ra \text{hocolim}_{I}cX$ is
an isomorphism in $\text{Ho}(\C)$ for any $X\in \C$. To prove the
proposition it remains to show that a weak equivalence of spaces
$f:K\ra L$ induces an isomorphism $ {\mathbf f}^{op}\otimes_{l}X:
{\mathbf K}^{op}\otimes_{l}X\ra {\mathbf L}^{op}\otimes_{l}X$ in
$\text{\rm Ho}(\C)$.

Let $\bf S$ be the collection of maps of spaces $f:K\ra L$ such
that, for any category  $(\C,W)$ admitting a left model
approximation and for any $X\in \C$, ${\mathbf
f}^{op}\otimes_{l}X$ is an isomorphism in $\text{Ho}(\C)$. We
claim that:
\begin{enumerate}
\item If $K$ and $L$ are contractible, then any $f:K\ra L$ belongs
to $\bf S$. \item Maps of the form $\Delta[n,k]\mono\Delta[n]$
belong to $\bf S$. \item $\bf S$ is closed under retracts. \item
$\bf S$ is closed under arbitrary sums. \item If an inclusion
$f:K\mono L$ belongs to $\bf S$, then, for any map $K\ra A$, so
does the inclusion $A\mono \text{colim}(A\la
K\stackrel{f}{\mono}L)$.
\item If a sequence of inclusions
$K_{i}\mono K_{i+1}$ ($i>0$) belong to $\bf S$, then so does
$K_{1}\mono \text{colim}_{i>0}K_{i}$.
\item A one-sided inverse of an element in $\bf S$ belongs to $\bf
S$.
\end{enumerate}
Note that properties (2)-(6) imply that any trivial cofibration of
spaces belongs to $\bf S$. And since any weak equivalence can be
factored as a trivial cofibration followed by a left inverse of a
trivial cofibration, the proposition clearly follows from
property~(7). Let us thus prove the claim.

Property (1) is a consequence of~\cite[Corollary
29.2]{MR2002k:55026} and (2) is a particular case of it.
Property~(3) is obvious and property (7) follows from a
``two-out-of-three argument''. Since the arguments for (4), (5),
and (6) are basically the same, we present the details of how to
prove (5) only.

Set $D$ to be the space $\text{colim}(A\la K \stackrel{f}{\mono}
L)$ and fix a left model approximation $l:\M\rightleftarrows\C:r$.
For an object $X\in \C$, choose a cofibrant replacement $QX$ in
$\text{Fun}^{b}\big(N(\mathbf D^{op}),\M\big)$ of the constant
diagram $crX$. Consider next the following push-out square  of
spaces:
\[
\xymatrix{ N(\mathbf K^{op})\rrmono^{N(\mathbf f^{op})}\dto && N(\mathbf L^{op})\dto^{h}\\
N(\mathbf A^{op})\rrmono && N(\mathbf D^{op})}
\]
Note that all the maps in this diagram are reduced and the
horizontal maps are cofibrations. It follows that $h^{\ast}QX$ is
cofibrant in $\text{Fun}^{b}\big(\mathbf N(\mathbf
L^{op}),\M\big)$ and consequently the morphism
$\text{colim}_{\mathbf N(\mathbf K^{op})}QX\ra
\text{colim}_{\mathbf N(\mathbf L^{op})}QX$ is a cofibration
(Proposition~\ref{prop ext}). Moreover it is an acyclic
cofibration as $f$ is assumed to be in ${\bf S}$. This shows that
$\text{colim}_{\mathbf N(\mathbf A^{op})}QX\ra
\text{colim}_{\mathbf N(\mathbf D^{op})}QX$ is also an acyclic
cofibration between cofibrant objects. The morphism
$l(\text{colim}_{\mathbf N(\mathbf A^{op})}QX)\ra
l(\text{colim}_{\mathbf N(\mathbf D^{op})}QX)$ is therefore a weak
equivalence in $\C$. Since it can be identified with ${\mathbf
A}^{op}\otimes_{l}X\ra {\mathbf D}^{op}\otimes_{l}X$, the claim is
proven.\hfill{\qed}

\section{An almost simplicial model structure}
\label{sec almost}
Let $\M$ be a model category and $A$ be a space. The category
$\text{Fun}^{b}\big({\mathbf N({\mathbf A})},\M\big)$ has both a
model structure (Theorem~\ref{thm modelmain}) and a simplicial
structure (Proposition~\ref{prop simplicial}). How compatible are
they?  To what extent is it a simplicial model category? How
related are the functors:
\[
\otimes_{l}:Spaces\times \text{Fun}^{b}\big({\mathbf N({\mathbf
A})},\M\big)\ra \text{Ho}\big(\text{Fun}^{b}({\mathbf N({\mathbf
A})},\M)\big)
\]
\[
\otimes: Spaces\times\text{Fun}^{b}\big({\mathbf N({\mathbf
A})},\M\big)\ra \text{Fun}^{b}\big({\mathbf N({\mathbf
A})},\M\big)?
\]

\begin{example}\label{ex negative}{\rm
Let us choose a simplex $\sigma:\Delta[0]\ra N(\mathbf A)$. By
definition it is given by an object in $\mathbf A$, i.e., a map
$\Delta[n]\ra A$. Let $F\in \text{Fun}^{b}\big({\mathbf N({\mathbf
A})},Spaces\big)$ be a cofibrant diagram such that
$F(\sigma)\not=\emptyset$. Note that:
\[
\text{lim}\big(\Delta[0]\xra{\sigma} N(\mathbf A)\xla{\pi}
N({\mathbf \Delta[1]}\boxtimes {\mathbf A})\big)
=\coprod_{\text{mor}(\mathbf \Delta[n],\mathbf
\Delta[1])}\Delta[0]
\]
Therefore $(\Delta[1]\otimes
F)(\sigma)=\pi^{k}\pi^{\ast}F=\coprod_{\text{mor}(\mathbf
\Delta[n],\mathbf \Delta[1])}F(\sigma)$ which is not weakly
equivalent to $F(\sigma)$.  It follows that the map $\Delta[1]\ra
\Delta[0]$ does not induce a weak equivalence of functors
$\Delta[1]\otimes F\ra \Delta[0]\otimes F=F$.

Thus the operation $-\otimes F$ does not preserve weak
equivalences, and consequently $\text{Fun}^{b}\big({\mathbf
N({\mathbf A})},Spaces\big)$ is not a simplicial model category.
Further, since the functor $-\otimes_{l}F$ preserves weak
equivalences (see Corollary~\ref{col hoinvcf}), the objects
$\Delta[1]\otimes F$ and ${\mathbf \Delta[1]}\otimes_{l} F$ are
not isomorphic in $\text{Ho}\big (\text{Fun}^{b}(\mathbf N(\mathbf
A),Spaces)\big)$}.
\end{example}

Despite this negative example, we are going to show that only half
of Quillen's axiom (SM7)~\cite{MR36:6480} does not hold for
$\text{Fun}^{b}\big(\mathbf N(\mathbf A),\M\big)$ to be a
simplicial model category. Moreover we are going to prove that,
for any homotopically constant  and cofibrant $F\in
\text{Fun}^{b}\big({\mathbf N({\mathbf A})},\M\big)$,
$\text{colim}_{\mathbf N({\mathbf A})}(K\otimes F)$ and
$\text{ocolim}_{\mathbf N({\mathbf A})}({\mathbf K}\otimes_{l} F)$
(see the remark after Theorem~\ref{thm hoinvcolim}) are isomorphic
in $\text{Ho}(\M)$, although $K\otimes F$ and  ${\mathbf
K}\otimes_{l} F$ might fail to be so in
$\text{Ho}\big(\text{Fun}^{b}({\mathbf N({\mathbf A})},\M)\big)$.

 Let $\phi:F\ra G$ be a
morphism in $\text{Fun}^{b}\big({\mathbf N({\mathbf A})},\M\big)$
and \mbox{$f:K \ra L$} be a map of spaces. Define:
\[
H:=\text{colim}(L\otimes F\xla{f\otimes F}  K\otimes
F\xra{K\otimes \phi}K\otimes G)
\]
Let $\psi:H\ra L\otimes G$ in  $\text{Fun}^{b}\big({\mathbf
N(\mathbf A)},\M\big)$ be induced by $f$, $\phi$, and the
universal property of the colimit. The following proposition shows
that half of Quillen's axiom (SM7)~\cite{MR36:6480} holds for the
simplicial structure introduced in Section~\ref{sec subdivision}.

\begin{proposition}\label{prop amodst}
If $f:K\mono L$ is an inclusion of spaces and $\phi:F\ra G$ is an
(acyclic) cofibration in $\text{Fun}^{b}\big({\mathbf N({\mathbf
A})},\M\big)$, then $\psi:H\ra L\otimes G$ is an (acyclic)
cofibration in $\text{Fun}^{b}\big({\mathbf N({\mathbf
A})},\M\big)$.
\end{proposition}

\begin{proof}
We need to show that, for a non-degenerate $\sigma:\Delta[n]\ra
N(\mathbf A)$, the morphism:
\[
i:\text{colim}\big(H(\sigma)\la\text{colim}_{\mathbold{\partial\Delta}[n]}H\ra
\text{colim}_{\mathbold{\partial\Delta}[n]}(L\otimes G)\big) \lra
(L\otimes G)(\sigma)
\]
is an (acyclic) cofibration in $\M$ (Theorem~\ref{thm modelmain}).
To identify both sides of this morphism consider the following
commutative diagram, where all the squares are pull-backs:
\[
\xymatrix{
\partial P \rmono\dmono & P\rto\dmono &  N({\mathbf K}\boxtimes {\mathbf A})\dmono\\
 \partial Q\rmono\dto  & Q\rto\dto^{p} &  N({\mathbf L}\boxtimes {\mathbf A})\dto\\
 \partial\Delta[n]\rmono  & \Delta[n]\rto^{\sigma}  & N(\mathbf A)
}
\]
All the vertical maps in this diagram are reduced, as they are
pull-backs of reduced maps (Proposition~\ref{prop simpl}).
Moreover, since $\partial Q \mono Q$ and $P \mono Q$ are
inclusions, the map $R:=\text{colim}(\partial Q \la
\partial P\ra P)\ra Q$ is a monomorphism. The fact that $p:Q\ra
\Delta[n]$ is reduced, implies that the only non-degenerate
simplices in $Q\setminus \partial Q$ are of dimension $n$ and $p$
sends them to the only non-degenerate simplex of dimension $n$ in
$\Delta[n]$. As $\sigma$ is non degenerate in $N(\mathbf A)$ and
$(Q\setminus R)\subset (Q\setminus \partial Q)$, it follows that
the composition $\sigma p$ takes non-degenerate simplices in
$Q\setminus R$ to non-degenerate simplices in $N({\mathbf A})$.

The spaces in this diagram can be used  to make the following
identifications:
\[\begin{array}{rclcrcl}
\text{colim}_{\mathbold{\partial\Delta}[n]}(L\otimes F) & = &
\text{colim}_{\mathbold{\partial Q}}F & \hspace{5pt} & (L\otimes
F)(\sigma) & = & \text{colim}_{\mathbf Q}F
\\
\text{colim}_{\mathbold{\partial\Delta}[n]}(L\otimes G) & =&
\text{colim}_{\mathbold{\partial Q}}G &  & (L\otimes G)(\sigma) &
= & \text{colim}_{\mathbf Q}G
\\
\text{colim}_{\mathbold{\partial\Delta}[n]}(K\otimes F) & =&
\text{colim}_{\mathbold{\partial P}}F &  & (K\otimes F)(\sigma) &
= &
\text{colim}_{\mathbf P}F\\
\text{colim}_{\mathbold{\partial\Delta}[n]}(K\otimes G) &
=&\text{colim}_{\mathbold{\partial P}}G &  & (K\otimes G)(\sigma)
& = & \text{colim}_{\mathbf P}G
\end{array}\]
Thus the left hand side of the morphism $i$ coincides with  the
colimit of the cube:
\[\xymatrix{
 & \text{colim}_{\mathbold{\partial P}}F\rrto\ddto|\hole\dlto
& & \text{colim}_{\mathbf P}F\dlto\ddto\\
\text{colim}_{\mathbold{\partial P}}G \rrto\ddto & & \text{colim}_{\mathbf P} G\\
 & \text{colim}_{\mathbold{\partial Q}}F\rrto\dlto & & \text{colim}_{\mathbf Q} F\\
 \text{colim}_{\mathbold{\partial Q}} G
}\] Note that:
\[
\text{colim}\big( \text{colim}_{\mathbold{\partial Q}}F\la
\text{colim}_{\mathbold{\partial P}}F\ra \text{colim}_{\mathbf
P}F\big)=\text{colim}_{\mathbf R}F
\]
\[
\text{colim}\big( \text{colim}_{\mathbold{\partial Q}}G\la
\text{colim}_{\mathbold{\partial P}}G\ra \text{colim}_{\mathbf
P}G\big)=\text{colim}_{\mathbf R}G
\]
Thomason's \cite[Theorem~26.8]{MR2002k:55026} permits in general
to decompose the (ho)colimit over a Grothendieck construction into
simpler (ho)colimits. In the particular case of a cube as above,
\cite[Example~38.2]{MR2002k:55026} shows that the colimit can be
also expressed as a pushout (see also
\cite[Claim~2.8]{MR93i:55015}):
\[
\text{colim}(\text{colim}_{\mathbf Q}F\la \text{colim}_{\mathbf
R}F\mono \text{colim}_{\mathbf R}G)
 \]
Moreover $i$ can be identified with the morphism induced by the
inclusion $R\subset Q$ and the  natural transformation $\phi:F\ra
G$:
\[
 \text{colim}(\text{colim}_{\mathbf Q}F\la \text{colim}_{\mathbf R}F\mono
 \text{colim}_{\mathbf R}G)\ra  \text{colim}_{\mathbf Q}G
\]
We can now apply Proposition~\ref{prop ext} to conclude that this
morphism is an (acyclic) cofibration in $\M$.
\end{proof}

Although in general the operation $K\mapsto K\otimes F$ is not
homotopy invariant (Example~\ref{ex negative}), it behaves well on
homotopically constant diagrams.

\begin{proposition}\label{prop colimhoinvcons}
Let $F\in \text{Fun}^{b}\big(\mathbf N(\mathbf A),\M\big)$ be
homotopically constant (\ref{pt functors}) and cofibrant. Then:
\begin{enumerate}
\item For any a weak equivalence of spaces $f:K\ra L$, the induced morphism
$\text{\rm colim}_{\mathbf N(\mathbf A)}(f\otimes F): \text{\rm
colim}_{\mathbf N(\mathbf A)}(K\otimes F)\ra \text{\rm
colim}_{\mathbf N(\mathbf A)}(L\otimes F) $ is a weak equivalence
in $\M$.
\item
There is an isomorphism $\text{\rm colim}_{\mathbf N(\mathbf A)}(K
\otimes F) \cong \text{\rm ocolim}_{\mathbf N(\mathbf A)}({\mathbf
K} \otimes_l F)$  in $\text{Ho}(\M)$, for any space $K$.
\end{enumerate}
\end{proposition}

\begin{proof}
Choose a weak equivalence  $F\ra cX$ in $
\text{Fun}^{b}\big({\mathbf N({\mathbf A})},\M\big)$. Let  $\pi$
be the subdivision of the projection $K\times A\ra A$ and note
that there is an isomorphism $\text{colim}_{\mathbf N({\mathbf
A})}(K\otimes F)\cong\text{colim}_{\mathbf N({\mathbf K}\boxtimes
{\mathbf A})}\pi^{\ast}F$. Since $\pi$ is reduced, the bounded
diagram $\pi^{\ast}F$ is cofibrant in $\text{Fun}^{b}\big({\mathbf
N({\mathbf K}\boxtimes {\mathbf A})},\M\big)$. It follows  that
there is an isomorphism $\text{colim}_{\mathbf N({\mathbf
A})}(K\otimes F)\cong ({\mathbf K}\boxtimes {\mathbf
A})\otimes_{l}X$ in $\text{Ho}(\M)$. By the same argument,   the
morphism $({\mathbf f}\boxtimes {\mathbf A})\otimes_{l}X:
({\mathbf K}\boxtimes {\mathbf A})\otimes_{l}X\ra ({\mathbf
L}\boxtimes {\mathbf A})\otimes_{l}X$ can be identified with
$\text{\rm colim}_{\mathbf N({\mathbf A})}(f\otimes F)$. Since
$f\times A$ is a weak equivalences, according to
Corollary~\ref{col hoinvcf}, $({\mathbf f}\boxtimes {\mathbf
A})\otimes_{l}X$ is an isomorphism in $\text{Ho}(\M)$ and so
$\text{\rm colim}_{\mathbf N({\mathbf A})}(f\otimes F)$ is a weak
equivalence in $\M$. This shows the first statement.

To prove the second statement consider the following sequence of
isomorphisms in $\text{Ho}(\M)$:
\begin{align*}
\text{colim}_{\mathbf N({\mathbf A})}(K\otimes F)& \cong ({\mathbf
K}\boxtimes {\mathbf A})\otimes_{l}X\stackrel{(a)}{\cong}
({\mathbf K}\times {\mathbf A})\otimes_{l}X\stackrel{(b)}{\cong}
({\mathbf A}\times {\mathbf K})\otimes_{l}X \\
& \stackrel{(c)}{\cong}{\mathbf A}\otimes_{l}({\mathbf
K}\otimes_{l}X)\cong\text{\rm ocolim}_{\mathbf N(\mathbf
A)}({\mathbf K} \otimes_l F)
\end{align*}
where $(a)$ and $(b)$ are isomorphisms because of
Corollary~\ref{col hoinvcf} and $(c)$ is an isomorphism because
of~\cite[Theorem 24.9]{MR2002k:55026}.
\end{proof}

How do these properties of the operation $\otimes$ in $\text{\rm
Fun}^{b}\big(\mathbf N(\mathbf A),\M\big)$ translate to the
mapping spaces introduced in Proposition~\ref{prop simplicial}?

\begin{corollary}\label{col mapfun}\hspace{1mm}
If $\phi$  is a weak equivalence between cofibrant objects and
$\psi$ is a weak equivalence between fibrant objects in $\text{\rm
Fun}^{b}\big(\mathbf N(\mathbf A),\M\big)$, then $\text{\rm
map}(\phi,\psi)$ is a weak equivalence of spaces.
\end{corollary}
\begin{proof}
By standard K. Brown lemma type of arguments \cite[Lemma
9.9]{MR96h:55014}, without loss of generality, we can assume in
addition that $\phi$ is an acyclic cofibration and $\psi$ is an
acyclic fibration.  In this case the corollary follows from
Proposition~\ref{prop amodst} and the usual adjointness argument.
\end{proof}

The mapping space functor satisfies also the following duality
property:

\begin{proposition}\label{prop basicduality}
If $H\in \text{\rm Fun}^{b}\big(\mathbf K, \text{\rm
Fun}^{b}(\mathbf N(\mathbf A),\M)\big)$ is a cofibrant functor
whose values are homotopically constant functors in $ \text{\rm
Fun}^{b}\big(\mathbf N(\mathbf A),\M\big)$ and $Z$ is fibrant in
$\M$, then $\text{\rm map}(H,cZ)$ is fibrant in $\text{\rm
Fun}^{b}(\mathbf K^{op},Spaces)$~\cite[Section
31.3]{MR2002k:55026}.
\end{proposition}

\begin{proof}
We need to show that, for a non-degenerate simplex
$\sigma:\Delta[n]\ra K$, we have a fibration of spaces
$\text{map}(H(\sigma),cZ)\ra \text{lim}_{(\mathbf \partial\mathbf
\Delta [n])^{op}}\text{map}(H,cZ)=
\text{map}(\text{colim}_{\mathbold{\partial \Delta[n]}}H,cZ)$. To
do that one needs to construct certain lifts. By adjointness this
problem reduces to showing that, for  a push-out square of the
form:
\[\xymatrix{
\Delta[k,l]\otimes (\text{colim}_{\mathbold{\partial\mathbf \Delta[n]}}H)\rto\dto_{\alpha} & \Delta[k,l]\otimes H(\sigma)\dto\\
\Delta[k]\otimes (\text{colim}_{\mathbold{\partial\mathbf
\Delta[n]}}H)\rto & P(\sigma) \rto &  \Delta[k] \otimes H(\sigma)}
\]
any natural transformation $P\ra cZ$ factors through $P \ra
\Delta[k]\otimes H(-)$. Since $cZ$ is a constant functor and $Z$
is fibrant, it is enough to prove that the morphism
$\text{colim}_{\mathbf N(\mathbf A)}P \ra \text{colim}_{\mathbf
N(\mathbf A)}\big(\Delta[k]\otimes H(\sigma)\big)$ is an acyclic
cofibration in $\M$. Cofibrancy is a consequence of
Propositions~\ref{prop amodst} and~\ref{prop
basicprophocolim}.(1). Using Proposition~\ref{prop
colimhoinvcons},  acyclicity  would follow if
$\text{colim}_{\mathbf N(\mathbf A)}\alpha$ is an acyclic
cofibration in $\M$. Since  $H$ is cofibrant, the same is true for
$\text{colim}_{\mathbold{\partial\mathbf \Delta[n]}}H$ and so
$\alpha$ is a cofibration by Proposition~\ref{prop amodst}. This
implies that $\text{colim}_{\mathbf N(\mathbf A)}\alpha$ is also a
cofibration. It remains to show that $\text{colim}_{\mathbf
N(\mathbf A)}\alpha$ is a weak equivalence.

Note that, for any $L$,  $L\otimes H$ is cofibrant in $\text{\rm
Fun}^{b}\big(\mathbf K, \text{\rm Fun}^{b}(\mathbf N(\mathbf
A),\M)\big)$ (another consequence of Proposition~\ref{prop
amodst}). The values of $H$ are cofibrant and homotopically
constant. We can therefore apply Proposition~\ref{prop
colimhoinvcons} to conclude that the natural transformation
$\text{colim}_{\mathbf N(\mathbf A)}(\Delta[k,l]\otimes H)\ra
\text{colim}_{\mathbf N(\mathbf A)}(\Delta[k]\otimes H)$ is a weak
equivalence between cofibrant objects in $\text{Fun}^{b}(\mathbf
K,\M)$ \cite[Proposition 24.2]{MR2002k:55026}. Consequently
$\text{colim}_{\mathbf K}\text{colim}_{\mathbf N(\mathbf
A)}(\Delta[k,l]\otimes H)\ra \text{colim}_{\mathbf K}
\text{colim}_{\mathbf N(\mathbf A)}(\Delta[k]\otimes H)$ is a weak
equivalence in $\M$. As this last morphism can be identified with
$\text{colim}_{\mathbf N(\mathbf A)}\alpha$, the proposition is
proven.
\end{proof}

This implies that the connected components of certain mapping
spaces coincide as expected with homotopy classes of morphisms.

\begin{corollary}\label{col kanhomclas}
If $Z$ is fibrant in $\M$ and $F$ is  homotopically constant and
cofibrant in $\text{\rm Fun}^b\big(\mathbf N(\mathbf A),\M\big)$,
then $\text{\rm map}(F,cZ)$ is Kan. Moreover there is a bijection,
natural in $F$ and $Z$, between the set of components
$\pi_{0}\text{\rm map}(F,cZ)$  and  the set of morphism $[F,cZ]$
in $\text{\rm Ho}(\text{\rm Fun}^b\big(\mathbf N(\mathbf
A),\M)\big)$.
\end{corollary}

\begin{proof}
The first part is a particular case of Proposition~\ref{prop
basicduality} when $K=\Delta[0]$. To prove the second part, we
will show that the map which assigns to the component of a
morphism $\phi:F\ra cZ$ in $\text{\rm map}(F,cZ)$  the homotopy
class of $\phi$ in $[F,cZ]$, is a bijection.

Since $F$ is cofibrant, the natural transformation
$\partial\Delta[1]\otimes F\ra \Delta[1]\otimes F$, induced by the
inclusion $\partial\Delta[1]\mono \Delta[1]$, is a cofibration.
Consider two natural transformations $\phi_{0}:F\ra cZ$ and
$\phi_{1}:F\ra cZ$, choose a cylinder object $H$ for $F$, and
consider the following commutative diagram:
\[\xymatrix{
\partial\Delta[1]\otimes F\dto_{\text{id}\coprod\text{id}}\drmono\rrto ^{\phi_{0}\coprod \phi_{1}} & & cZ\\
F&\Delta[1]\otimes F\dmono\\
& H\ulto_{\sim} }
\]
We claim that there is a natural transformation $\Delta[1]\otimes
F\ra cZ$ which makes the above diagram commutative if and only if
there is a natural transformation $H\ra cZ$ with the same
property. Since the target $cZ$ is a constant functor any such map
$\Delta[1]\otimes F\ra cZ$ is induced by a morphism
$\text{colim}_{\mathbf N(\mathbf A)}(\Delta[1]\otimes F)\ra Z$ in
$\M$. From Proposition~\ref{prop colimhoinvcons} and
Proposition~\ref{prop basicprophocolim}(1) we infer that the
morphisms $\text{colim}_{\mathbf N(\mathbf A)}(\Delta[1]\otimes
F)\ra \text{colim}_{\mathbf N(\mathbf A)}F$ and
$\text{colim}_{\mathbf N(\mathbf A)}H\ra \text{colim}_{\mathbf
N(\mathbf A)}F$ are weak equivalences. It follows that the
morphism $\text{colim}_{\mathbf N(\mathbf A)}(\Delta[1]\otimes
F)\ra \text{colim}_{\mathbf N(\mathbf A)}H$ is an acyclic
cofibration in $\M$ and our claim is proven, as $Z$ is assumed to
be fibrant. This argument shows that $\phi_{0}$ and $\phi_{1}$ are
left homotopic if and only if they are in the same component of
$\text{\rm map}(F,cZ)$.
\end{proof}

\section{Resolutions and mapping spaces}
\label{sec resolution}
A standard way of studying derived or homotopical properties of
various operations is by considering resolutions. This resolution
principle is playing a significant role in our approach to mapping
spaces too. In our setting it is stated as a model approximation.

\begin{proposition}\label{prop resolution}
Let $\M$ be a model category. If  $A$ is a contractible space,
then the pair of adjoint functors $\text{\rm colim}_{\mathbf
N(\mathbf A)}:\text{\rm Fun}^{b}\big(\mathbf N(\mathbf
A),\M\big)\rightleftarrows \M:c$ is a  left model approximation of
$\M$.
\end{proposition}

\begin{proof}
Conditions (1), (2), and (3) of Definition~\ref{def modelapprox}
are obviously satisfied. Condition~(4) is a consequence of
Theorem~\ref{thm hoinvcolim}. To prove (5) consider a weak
equivalence $F\ra cX$ in $\text{\rm Fun}^{b}\big(\mathbf N(\mathbf
A),\M\big)$ with a cofibrant domain. According to Theorem~\ref{thm
hocolimconst}, $\text{colim}_{\mathbf N(\mathbf A)}F$ is
isomorphic to ${\mathbf A}\otimes_{l}X$ in $\text{Ho}(\M)$. Since
$A$ is contractible, the functor from the simplex category of $A$
to the category with only one object and one morphism is a weak
equivalence.  Thus, by Proposition~\ref{prop hoinvop}, $\mathbf
A\otimes_{l}X\ra X$ is also an isomorphism in $\text{Ho}(\M)$ and
so $\text{colim}_{\mathbf N(\mathbf A)}F\ra X$ is a weak
equivalence in $\M$.
\end{proof}

Taking the subdivision of $A$ is a crucial assumption in the
statement of Proposition~\ref{prop resolution}. It is not true in
general that if $A$ is contractible, then the pair of adjoint
functors $\text{colim}_{\mathbf A}: \text{Fun}^{b}(\mathbf A,\M)
\rightleftarrows \M:c$ forms a left model approximation. Take $\M$
to be $Spaces$ and $A$ to be the space defined in  Example~\ref{ex
nonhoinv}. Consider a weak equivalence $F\ra c\Delta[0]$ in
$\text{Fun}^{b}(\mathbf A,Spaces)$ with a cofibrant domain. The
space $\text{colim}_{\mathbf A} F$ is weakly equivalent to
$\text{hocolim}_{\mathbf A[S^{-1}]} \Delta[0] \simeq S^{2}$, which
is not weakly contractible. Thus the requirement (5) of
Definition~\ref{def modelapprox} is not satisfied.

How to construct mapping spaces in a category $\C$ that admits a
left model approximation? Here is the recipe:
\begin{itemize}
\item choose a left model approximation $l:\M\rightleftarrows \C:r$;
\item choose a contractible space $A$;
\item choose a cofibrant replacement $Q$ in $\text{\rm Fun}^{b}\big(\mathbf N(\mathbf A),\M\big)$
and a fibrant replacement $R$ in $\M$;
\item  define a functor $\text{map}:\C^{op}\times \C\ra Spaces$ by the formula:
\[\text{map}(X,Y):=\text{map}\big(Qcr(X),cRr(Y)\big)\]
where the mapping space on the right comes from the simplicial
structure on $\text{\rm Fun}^{b}\big({\mathbf N({\mathbf
A})},\M\big)$ introduced in Proposition~\ref{prop simplicial}.
\end{itemize}

These mapping spaces have the following properties:
\begin{proposition}\label{prop basicmaplma}\hspace{1mm}
\begin{enumerate}
\item For any $X$ and $Y$ in $\C$,  $\text{\rm map}(X,Y)$ is Kan.
\item There is a bijection, natural in $X$ and $Y$ in $\C$, between the set
of connected components $\pi_{0}\text{\rm map}(X,Y)$ and the set
of morphism $[X,Y]$ in $\text{\rm Ho}(\C)$.
\item If $\alpha:X_{0}\ra X_{1}$ and $\beta:Y_{0}\ra Y_{1}$ are weak
equivalences in $\C$, then $\text{\rm map}(\alpha,\beta):\text{\rm
map}(X_{1},Y_{0})\ra \text{\rm map}(X_{0},Y_{1})$ is a weak
equivalence of spaces.
\end{enumerate}
\end{proposition}
\begin{proof}
Statement (1) is a direct consequence of the first part of
Corollary~\ref{col kanhomclas}.

According to the second part of this corollary there is a natural
bijection between the set of connected components
$\pi_{0}\text{\rm map}(X,Y)$ and the set of morphism
$[cr(X),cr(Y)]$ in $\text{Ho}\big(\text{\rm Fun}^{b}(\mathbf
N(\mathbf A),\M)\big)$. Since $A$ is contractible,
Proposition~\ref{prop resolution} implies that the pair
$l\text{colim}_{\mathbf N(\mathbf A)}:\text{\rm
Fun}^{b}\big(\mathbf N(\mathbf A),\M\big)\rightleftarrows \C:cr$
forms a left model approximation of $\C$. It follows that there is
a bijection between the set of morphisms $[X,Y]$ in
$\text{Ho}(\C)$ and $[cr(X),cr(Y)]$ in $\text{Ho}\big(\text{\rm
Fun}^{b}(\mathbf N(\mathbf A),\M)\big)$ \cite[Proposition
5.5]{MR2002k:55026}, which shows (2).

Statement (3) is a direct consequence of Corollary~\ref{col
mapfun}.
\end{proof}

\begin{remark}
\label{rem associative}
{\rm Assume that, in the above recipe for constructing mapping
spaces, we chose a model category $\mathcal M$ which has a
\emph{functorial} cofibrant-fibrant replacement. This induces in
turn a functorial cofibrant-fibrant replacement in $\text{\rm
Fun}^{b}\big(\mathbf N(\mathbf A),\M\big)$. For any $X$ in $\C$
define $\widetilde X$ to be this replacement for the constant
functor $cr(X)$. Such choice leads not only to a functor
$\text{map}(X,Y):=\text{map}\big(\widetilde X, \widetilde Y
\big)$, but also to a \emph{strictly associative composition}
\[
\text{map}(X,Y) \times \text{map}(Y,Z) \longrightarrow
\text{map}(X,Z)
\]
which comes from the simplicial enrichment introduced in
Proposition~\ref{prop simplicial}.}
\end{remark}

Part (3) of Proposition~\ref{prop basicmaplma} implies that the
mapping space functor induces a functor on the level of homotopy
categories, denoted by the same symbol:
\[
\text{map}:\text{Ho}(\C)^{op}\times
\text{Ho}(\C)\ra\text{Ho}(Spaces)
\]

How dependent is this functor on the choices we have made: the
model approximation $l:\M\rightleftarrows:\C:r$, the contractible
space $A$, and the cofibrant and fibrant replacements? It can be
characterized by the following universal property:

\begin{proposition}
\label{prop mapasadjoint}
Let $X\in \C$. The functor $\text{\rm map}(X,-):\text{\rm
Ho}(\C)\ra \text{\rm Ho}(Spaces)$ is right adjoint to
$-\otimes_{l} X:\text{\rm Ho}(Spaces)\ra \text{\rm Ho}(\C)$. In
particular it does not depend on  the choices of the model
approximation $l:\M\rightleftarrows\C:r$, the contractible space
$A$, and  the cofibrant and fibrant replacements.
\end{proposition}

\begin{proof}
Consider the following sequence of adjoint functors:
\[\xymatrix{
Spaces \ar@<0.7ex>[rr]^(0.4){-\otimes QcrX} & &
\text{Fun}^{b}\big(\mathbf N(\mathbf A),\M\big)
\ar@<0.7ex>[ll]^(0.6){\text{map}(QcrX,-)}\ar@<0.7ex>[rr]^(0.65){\text{colim}_{{\mathbf
N({\mathbf A})}}} & &\M
\ar@<0.7ex>[ll]^(0.35){c}\ar@<0.7ex>[r]^(0.6){l} &
\C\ar@<0.7ex>[l]^(0.38){r} }
\]
According to Propositions~\ref{prop basicprophocolim},~\ref{prop
amodst}, and~\ref{prop colimhoinvcons},  the functor
$\text{colim}_{\mathbf N(\mathbf A)}(-\otimes QcrX)$ preserves
cofibrations and acyclic cofibrations. Thus, by~\cite[Theorem
9.7]{MR96h:55014}, the total derived functors also form  adjoint
pairs (see also Lemma~\ref{lem adjmodapprox}):
\[
\xymatrix{ \text{Ho}(Spaces)
\ar@<0.7ex>[rrrr]^(0.53){L(\text{colim}_{\mathbf N(\mathbf
A)}(-\otimes QcrX))} & & & &
\text{Ho}(\M)\ar@<0.7ex>[llll]^(0.45){R(\text{map}(QcrX,c-))}
\ar@<0.7ex>[r]^(0.5){L(l)} & \text{Ho}(\C)\ar@<0.7ex>[l]^(0.5){r}}
\]
Note that the composition of $L(\text{colim}_{\mathbf N(\mathbf
A)}(-\otimes QcrX))$ and $L(l)$ can be identified with $K\mapsto
(\mathbf K\boxtimes \mathbf A) \otimes_{l}X$ (see the proof of
Proposition~\ref{prop colimhoinvcons}), and the composition of $r$
and $R(\text{map}(QcrX,c-))$ with $\text{\rm map}(X,-):\text{\rm
Ho}(\C)\ra \text{\rm Ho}(Spaces)$. Since $A$ is assumed to be
contractible, the morphism $(\mathbf K\boxtimes \mathbf
A)\otimes_{l}X\ra {\mathbf K}\otimes_{l}X$, induced by the
projection $\pi:K\times A\ra K$, is an isomorphism in
$\text{Ho}(\C)$ (Corollary~\ref{col hoinvcf}).
\end{proof}

\begin{corollary}
Assume that $(\C,W)$ admits a left model approximation. The
functor $\otimes_{l}:\text{\rm Ho}(Spaces)\times \text{\rm
Ho}(\C)\ra \text{\rm Ho}(\C)$ is then a homotopy left action
(\ref{pt hosimpstruct}). \hfill{\qed}
\end{corollary}

How natural is this homotopy left action? Let $(\C,W)$ and
$(\D,V)$ admit left model approximations and $\Phi:\C\ra \D$ be a
functor that sends weak equivalences in $\C$ to weak equivalences
in $\D$. The universal property of the total left derived functor
$\text{hocolim}_{I}$ yields a morphism $\text{hocolim}_{I}\Phi
F\ra \Phi(\text{hocolim}_{I}F)$ in $\text{Ho}(\D)$ natural in
$F:I\ra \C$ (see~\ref{pt nathocolim}). For $X\in \C$ and $K\in
Spaces$, define:
\[{\mathbf K}\otimes_{l}\Phi (X)\ra
\Phi({\mathbf K}\otimes_{l}X)\] to be the morphism
$\text{hocolim}_{\mathbf K}c\Phi(X)\ra
\Phi(\text{hocolim}_{\mathbf K}cX)$ in $\text{Ho}(\D)$.

\begin{corollary}
Assume that $(\C,W)$ and $(\D,V)$ admit left model approximations
and $\Phi:\C\ra \D$ is a functor that sends weak equivalences in
$\C$ to weak equivalences in $\D$. Then the induced functor
$\Phi:\text{\rm Ho}(\C)\ra \text{\rm Ho}(\D)$ together with the
morphisms ${\mathbf K}\otimes_{l}\Phi (X)\ra \Phi({\mathbf
K}\otimes_{l}X)$ is a simplicial functor (\ref{pt
homsimpfun}).\hfill{\qed}
\end{corollary}

We prove finally the classical property relating mapping spaces
out of a homotopy colimit with the corresponding homotopy limit of
mapping spaces, see \cite[XII Proposition~4.1]{MR51:1825} and
\cite[Theorem~19.4.4]{MR1944041}.

\begin{proposition}
\label{prop maphocolim}
Assume that $(\C, W)$ admits a left model approximation. Consider
an object $Y$ in $\C$ and $F:I\ra \C$ a functor. Then there is an
isomorphism $\text{\rm map}(\text{\rm hocolim}_{I}F,Y)\cong
\text{\rm holim}_{I^{op}}\text{\rm map}(F,Z)$ in $\text{\rm
Ho}(Spaces)$ natural in  $Y$ and $F$.
\end{proposition}

\begin{proof}
Let us choose a left model approximation $l:\M\rightleftarrows
\C:r$ and a contractible space $A$. Apply then
Proposition~\ref{prop basicduality} to a fibrant replacement
$Z\in\M$ of $rY$ and a cofibrant replacement $H\in
\text{Fun}^{b}\big(\mathbf N(I),\text{Fun}^{b}(\mathbf N(\mathbf
A),\M)\big)$ of the bounded diagram $\mathbf N(I)\xra{\epsilon}
I\xra{F} \C\xra{r} \M\xra{c} \text{Fun}^{b}(\mathbf N(\mathbf
A),\M)$.
\end{proof}

\section{Right model approximations and duality}
\label{sec right}
In this section we study categories admitting both a left and a
right model approximations. These structures lead to homotopy left
and right actions. How related are these actions? We prove below
in Theorem~\ref{prop homotopyadjoint} that the corresponding
simplicial structures on $\text{\rm Ho}(\C)$ are isomorphic (and
in particular the mapping spaces induced by the homotopy left and
right actions are isomorphic). We need first some dual statements
to those we have proven so far for left model approximations. They
can be restated by simply taking the opposite categories.

\begin{point}{\rm
One can define a simplicial structure on
$\text{Fun}^{b}\big(\mathbf N(\mathbf A)^{op},\C\big)$ by
identifying it with $\text{Fun}^{b}\big(\mathbf N(\mathbf
A),\C^{op}\big)$ and applying Proposition~\ref{prop simplicial}.
Explicitly, this simplicial structure is given by a functor:
\[
\text{hom}: Spaces^{op}\times \text{Fun}^{b}\big(\mathbf N(\mathbf
A)^{op},\C\big)\ra \text{Fun}^{b}\big(\mathbf N(\mathbf
A)^{op},\C\big)
\]
described as follows. Denote by
$\pi_{k}:\text{Fun}^{b}\big(\mathbf N(\mathbf K\boxtimes \mathbf
A)^{op},\C\big)\ra \text{Fun}^{b}\big( \mathbf N(\mathbf
A)^{op},\C\big)$ the right adjoint of
$\pi^{\ast}:\text{Fun}^{b}\big(\mathbf N(\mathbf
A)^{op},\C\big)\ra \text{Fun}^{b}\big(\mathbf N(\mathbf K\boxtimes
\mathbf A)^{op},\C\big)$, where as before $\pi:\mathbf N(\mathbf
K\boxtimes \mathbf A)\ra \mathbf N(\mathbf A)$ stands for the
subdivision of the projection onto the second factor. Define then
$\text{hom}(K,F):=\pi_{k}\pi^{\ast}F $. }\end{point}

\begin{point}{\rm
Assume that $(\C,W)$  admits a right model approximation. Then it
can be right localized. For any small category $I$, $\text{\rm
holim}_{I}:\text{\rm Ho}\big(\text{\rm Fun}(I,\C)\big)\ra
\text{\rm Ho}(\C)$ is right adjoint to $c:\text{\rm Ho}(\C)\ra
\text{\rm Ho}\big(\text{\rm Fun}(I,\C)\big)$. The assignment
$(I,X)\mapsto \text{\rm holim}_{I}cX$ defines a functor $\text{\rm
rhom}:Cat^{op}\times \C\ra \text{\rm Ho}(\C)$. This functor  is
homotopy invariant: if $f:I\ra J$ is a weak equivalence of small
categories, then $\text{\rm rhom}(f,X)$ is an isomorphism in
$\text{Ho}(\C)$ for any $X\in \C$ (compare with
Proposition~\ref{prop hoinvop}). }\end{point}

\begin{point}\label{pt rightmap}{\rm
If $f:K\ra L$ is a weak  equivalence of spaces, then the induced
morphism $\text{rhom}(\mathbf f,X):\text{rhom}(\mathbf L,X)\ra
\text{rhom}(\mathbf K,X)$ is an isomorphism in $\text{Ho}(\C)$,
for any $X\in \C$. The functor $\text{\rm rhom}:Spaces^{op}\times
\C\ra \text{Ho}(\C)$ is therefore a composition of the
localization $Spaces^{op}\times \C\ra \text{\rm Ho}(Spaces)^{op}
\times\text{\rm Ho}(\C)$ and a functor denoted by the same symbol
$\text{\rm rhom}:\text{\rm Ho}(Spaces)^{op}\times\text{\rm
Ho}(\C)\ra\text{\rm Ho}(\C) $. This functor is a homotopy right
action (\ref{pt homsimpfun}), i.e., it has the following
properties, dual to those in Corollary~\ref{lem prophococons} and
Proposition~\ref{prop mapasadjoint}:
\begin{enumerate}
\item[(1d)]
 there are isomorphisms $X\ra \text{\rm rhom}(\mathbf \Delta[0],X)$
in $\text{\rm Ho}(\C)$  natural in $X$;
\item[(2d)]
 there are coherent isomorphisms $\text{\rm rhom}({\mathbf K}\boxtimes {\mathbf L},X)\cong
\text{\rm rhom}\big(\mathbf K,\text{\rm rhom}(\mathbf L,X)\big)$
in $\text{\rm Ho}(\C)$ natural in $X\in \text{\rm Ho}(\C)$ and
$K,L\in \text{\rm Ho}(Spaces)^{op}$;
\item[(3d)] for any $X\in \text{\rm Ho}(\C)$, $\text{rhom}(-, X):
\text{\rm Ho}(Spaces)^{op}\ra \text{\rm Ho}(\C)$ has a left
adjoint $\text{map}(-,X):\text{\rm Ho}(\C)\ra \text{\rm
Ho}(Spaces)^{op}$.
\end{enumerate}
}\end{point}

\begin{point}{\rm
Let $(\C,W)$ and $(\D,V)$ be categories that admit right model
approximations and $\Phi:\C\ra \D$ be a functor that sends weak
equivalences in $\C$ to weak equivalences in $\D$.  This  and the
universal property of $\text{holim}$, as a total right derived
functor, implies that, for any $F:I\ra \C$, there is a morphism
$\Phi(\text{holim}_{I}F)\ra \text{holim}_{I}\Phi F$ in
$\text{Ho}(\D)$ natural in $F$ (see~\ref{pt nathocolim}). For
$X\in \C$ and $K\in Spaces$, define $\Phi\big(\text{rhom}(\mathbf
K,X)\big)\ra \text{rhom}(\mathbf K,\Phi F)$ to be the morphism
$\Phi(\text{holim}_{\mathbf K}cX)\ra \text{holim}_{\mathbf K}c\Phi
X$ in $\text{Ho}(\D)$. In this way the functor
$\Phi:\text{Ho}(\C)\ra \text{Ho}(\D)$ together with the morphism
$\Phi\big(\text{rhom}(\mathbf K,X)\big)\ra \text{rhom}(\mathbf
K,\Phi F)$ is a simplicial functor. }\end{point}

\begin{point}{\rm
The functor $\text{map}(-,X):\text{\rm Ho}(\C)\ra \text{\rm
Ho}(Spaces)^{op}$, defined in~\ref{pt rightmap}, can be rigidified
as follows:
\begin{itemize}
\item choose a right model approximation $l:\C\rightleftarrows \M:r$;
\item choose a contractible space $A$;
\item choose a fibrant replacement $R$ in $\text{\rm Fun}^{b}\big(\mathbf N(\mathbf A)^{op},\M\big)$
and a cofibrant replacement $Q$ in $\M$;
\item define a functor $\text{map}:\C^{op}\times \C\ra Spaces$ by the formula:
\[\text{map}(X,Y):=\text{map}\big(cQl(X),Rcr(Y)\big)\]
where the mapping space on the right comes from the simplicial
structure on $\text{\rm Fun}^{b}\big(\mathbf N(\mathbf
A)^{op},\M\big)$.
\end{itemize}
The functor $\text{map}:\C^{op}\times \C\ra Spaces$, defined by
the above procedure, is homotopy invariant and its values are Kan
spaces. The induced functor on the level of homotopy categories is
naturally isomorphic to the functor described in~\ref{pt rightmap}
and therefore it does not depend on the choice of the right model
approximation $l:\C\rightleftarrows \M:r$, the contractible space
$A$, and the fibrant and cofibrant replacements. }\end{point}

\begin{point}{\rm
The mapping space into a homotopy limit satisfies a property
analogous to Proposition~\ref{prop maphocolim}. Let $Y\in\C$ and
$F:I\ra \C$ be a functor. Then there is an isomorphism $\text{\rm
map}(Y, \text{\rm holim}_{I}F)\cong \text{\rm holim}_{I}\text{\rm
map}(Z,F)$ in $\text{\rm Ho}(Spaces)$ natural in $Y$ and $F$.
}\end{point}

We are now ready to compare the homotopy left and right actions
arising from left and right model approximations.

\begin{theorem}
\label{prop homotopyadjoint}
Assume that $(\C,W)$ admits both a left and a right model
approximations. Then, for any space $K$, the functor ${\mathbf
K}\otimes_{l} -:\text{\rm Ho}(\C)\ra \text{\rm Ho}(\C)$ is left
adjoint to $\text{\rm rhom}(\mathbf K,-):\text{\rm Ho}(\C)\ra
\text{\rm Ho}(\C)$. The weak simplicial structures on $\text{\rm
Ho}(\C)$ (\ref{pt hosimpstruct}), induced by the homotopy actions
$\otimes_{l}$  and $\text{\rm rhom}$, are isomorphic. In
particular the induced mapping spaces are also isomorphic.
\end{theorem}

\begin{proof}
According to Corollary~\ref{col althocolim}, we have the following
adjoint pairs:
\[
\xymatrix{ \text{Ho}(\C) \ar@<0.7ex>[rr]^(0.4){c}  &
&\text{Ho}\big(\text{Fun}(\mathbf
K,\C)\big)\ar@<0.7ex>[ll]^(0.58){\text{holim}_{\mathbf K}}
\ar@<0.7ex>[rr]^(0.6){\text{hocolim}_{\mathbf K}}  & &
\text{Ho}(\C)\ar@<0.7ex>[ll]^(0.4){c}}
\]
Since, by definition, $\text{hocolim}_{\mathbf K}cX=\mathbf
K\otimes_{l}X$ and $\text{holim}_{\mathbf K}cX=\text{rhom}(\mathbf
K,X)$ the theorem is proven.
\end{proof}

\section{Monoidal model categories}
\label{sec monoidal}

Recall that a monoidal model category is a model category $\M$
together with a ``compatible'' monoidal structure $\M \otimes \M
\rightarrow \M$, see \cite[Definition~4.2.6]{99h:55031}. Hovey
asks whether the homotopy category of such a monoidal model
category always forms a central algebra over $\text{\rm
Ho}(Spaces)$, i.e. if there is a monoidal functor $i: \text{\rm
Ho} (Spaces) \rightarrow \text{\rm Ho}(\M)$ and natural
isomorphisms $t: i(K) \otimes X \rightarrow X \otimes i(K)$
satisfying certain coherence conditions introduced in
\cite[Definition~4.1.10]{99h:55031}. We show this is indeed so.

\begin{proposition}
\label{prop monoidal}
Let $\M$ be a monoidal model category. Then the functor
$\otimes_{l}:\text{\rm Ho}(Spaces)\times \text{\rm
Ho}(\M)\ra\text{\rm Ho}(\M)$ describes a central algebra over
$\text{\rm Ho}(Spaces)$.
\end{proposition}
\begin{proof}
Let  us denote by $\otimes_{\M}:\M\times\M\ra \M$  the monoidal
structure in $\M$ and by $S \in \M$ the unit. By definition, for
any cofibrant $X\in \M$, the operations \mbox{$X\otimes_{\M}-$}
and $-\otimes_{\M}X$ preserve cofibrations. Consequently, for any
such $X$, if $F\in \text{Fun}^{b}(A,\M)$ is cofibrant, then so are
the functors $X\otimes_{\M}F$ and $F\otimes_{\M}X$.

Let $F:I\ra \M$ be a functor and  $QF\in
\text{Fun}^{b}\big(\mathbf N(I),\M\big)$ be a cofibrant
replacement of the composition of $\epsilon:\mathbf N(I)\ra I$ and
$F$. The above remark implies that, for a cofibrant $X$,
$QF\otimes_{\M}X$ and $X\otimes_{\M}QF$ are cofibrant replacements
of $F\epsilon\otimes_{\M}X$ and $X\otimes_{\M}F\epsilon$
respectively. Since $\otimes_{\M}$ commutes with colimits we have
isomorphisms in $\M$:
\[\text{colim}_{\mathbf N(I)}(QF\otimes_{\M}X)\cong (\text{colim}_{\mathbf N(I)}QF)\otimes_{\M}X\]
\[\text{colim}_{\mathbf N(I)}(X\otimes_{\M}QF)\cong X\otimes_{\M}(\text{colim}_{\mathbf N(I)}QF)\]
The corresponding isomorphisms in $\text{Ho}(\M)$ are easily seen
to be natural in $F:I\ra \M$ and $X\in\text{Ho}(\M)$:
\[\text{hocolim}_{I}(F\otimes_{\M}X)\cong (\text{hocolim}_{I}F)\otimes_{\M}X\]
\[\text{hocolim}_{I}(X\otimes_{\M}F)\cong X\otimes_{\M}(\text{hocolim}_{I}F)\]

It is now straightforward to check that the following isomorphism
$t$:
\[\xymatrix{
({\mathbf K}\otimes_{l}S)\otimes_{\M}X=(\text{hocolim}_{\mathbf
K}cS)\otimes_{\M}X &
\text{hocolim}_{\mathbf K}c(S\otimes _{\M}X)\lto\dto\\
 & \text{hocolim}_{\mathbf K}cX\\
X\otimes_{\M}({\mathbf
K}\otimes_{l}S)=X\otimes_{\M}(\text{hocolim}_{\mathbf K}cS)&
\text{hocolim}_{\mathbf K}c(X\otimes_{\M}S)\lto\uto }
\]
satisfies the coherence conditions for $\text{Ho}(\M)$ to be
central, where the monoidal functor $i$ is defined by $i(K) =
{\mathbf K}\otimes_{l}S$.
\end{proof}

%
%
%
%

\end{document}